\documentclass[11pt]{article}
\usepackage[all]{xy}
\usepackage{amsmath}
\usepackage{amsthm}
\usepackage{amscd}
\input amssym.tex

\input{epsf.tex}

\newcommand{\centeredepsfbox}[1]{\centerline{\epsfbox{#1}}}
\newcommand{\nb}[1]{#1\nobreakdash-}

\newcommand{\textmatrix}[4]{\bigl( \begin{smallmatrix} #1 & #2 \\ #3 & #4 
\end{smallmatrix} \bigr)}

\theoremstyle{definition}
\newtheorem*{Definition}{Definition}
\newtheorem*{Remark}{Remark}
\newtheorem*{Example}{Example}

\theoremstyle{plain}
\newtheorem{theorem}{Theorem}[section]
\newtheorem{proposition}[theorem]{Proposition}
\newtheorem{lemma}[theorem]{Lemma}
\newtheorem{corollary}[theorem]{Corollary}

\newcounter{remarks}
\newenvironment{Remarks}%
{\paragraph*{Remarks}\smallskip
     \begin{list}{\arabic{remarks}. }{\usecounter{remarks}%
          \setlength{\leftmargin}{0in}%
          \setlength{\rightmargin}{0in}%
          \setlength{\labelsep}{0pt}%
          \setlength{\labelwidth}{0pt}%
          \setlength{\listparindent}{0pt}%
     }
}
{
\end{list}
}

\DeclareMathOperator{\Ends}{Ends}

\DeclareMathOperator{\Aff}{Aff}

\DeclareMathOperator{\vcd}{vcd}
\DeclareMathOperator{\cd}{cd}
\DeclareMathOperator{\Transl}{Transl}
\DeclareMathOperator{\Stretch}{Stretch}

\DeclareMathOperator{\Length}{Length}

\DeclareMathOperator{\SL}{SL}

\DeclareMathOperator{\Homeo}{Homeo}
\DeclareMathOperator{\QIMap}{QIMap}
\DeclareMathOperator{\QC}{QC}
\DeclareMathOperator{\BS}{BS}
\DeclareMathOperator{\Isom}{Isom}

\DeclareMathOperator{\Sim}{Sim}
\DeclareMathOperator{\Aut}{Aut}
\DeclareMathOperator{\Stab}{Stab}

\DeclareMathOperator\Bilip{Bilip}
\DeclareMathOperator\Max{Max}

\DeclareMathOperator\radius{radius}
\DeclareMathOperator\Conf{Conf}

\newcommand\R{{\mathbf R}}
\newcommand\reals{\R}
\newcommand\Q{{\mathbf Q}}

\renewcommand\H{{\mathbf H}}
\newcommand\hyp{\H}
\newcommand\C{{\mathbf C}}
\newcommand\complex{\C}
\newcommand\Z{{\mathbf Z}}
\newcommand\integers{\Z}
\newcommand\naturals{{\mathbf N}}

\newcommand\inject{\hookrightarrow}
\newcommand\homeo{\approx}
\newcommand\isomorphic{\approx}
\newcommand\Sum{\sum}
\newcommand\infinity{\infty}
\newcommand\bndry{\partial}
\newcommand{\bdy}{\bndry}
\newcommand{\from}{\colon}
\def\composed{\circ}
\newcommand\suchthat{\bigm|}
\newcommand\inverse{{-1}}
\def\tensor{\otimes}
\newcommand\union{\cup}
\newcommand\Intersection{\cap}

\newcommand\mapright[1]{\xrightarrow{#1}}
\newcommand\absvalue[1]{\left| #1 \right|}

\newcommand\Id{\text{Id}}
\newcommand\A{\mathcal A}
\newcommand\intersect{\cap}

\newcommand\Svarc{\v{S}varc}

\newcommand\Solv{\text{Sol}}
\newcommand\restrict{\bigm|}
\newcommand\subgroup{<}
\newcommand\semidirect{\rtimes}
\newcommand\low{\ell}

\DeclareMathOperator\QI{QI}
\newcommand\cross{\times}

\newcommand\Hequiv{\sim}

\newcommand\Haus{H}

\begin{document}

\setcounter{page}{0}

\title{Quasi-isometric rigidity for the solvable Baumslag-Solitar
  groups, II.} 
\author{
Benson Farb and Lee Mosher
}
\date{Revised, December 10, 1997}

\thanks{The first author is supported in part by NSF grant DMS 9604640,
the second author by NSF grant DMS 9504946.}

\maketitle

\section*{Introduction}

Gromov's Polynomial Growth Theorem \cite{Gromov:PolynomialGrowth}
characterizes the class of virtually nilpotent groups by their
asymptotic geometry. Since Gromov's theorem it has been a major
open question (see, e.g.\ \cite{GhysHarpe:afterGromov}) to find an
appropriate generalization for solvable groups. This paper gives
the first step in that direction.

One fundamental class of examples of finitely-generated solvable
groups which are not virtually nilpotent are the \emph{solvable
Baumslag-Solitar groups}
$$ \BS(1,n) = \left< a,b \suchthat b a b^\inverse = a^n \right>
$$
where $n \ge 2$.  Our main theorem characterizes the group
$\BS(1,n)$ among all finitely-generated groups by its asymptotic
geometry.

\newtheorem*{QIRigidityTheorem}{Theorem A (Quasi-isometric rigidity)}
\begin{QIRigidityTheorem}
Let $G$ be any finitely generated group. If $G$ is quasi-isometric to
$\BS(1,n)$ for some $n \ge 2$, then there is a short exact sequence
$$1 \to N \to G \to \Gamma \to 1
$$
where $N$ is finite and $\Gamma$ is abstractly commensurable to $\BS(1,n)$.
\end{QIRigidityTheorem}

In fact we will describe the precise class of quotient groups
$\Gamma$ which can arise, and will classify all torsion-free $G$;
see section \ref{section:proof} in the outline below.

Theorem A complements the main theorem of \cite{FarbMosher:BSOne},
where it is shown that $\BS(1,n)$ is quasi-isometric to $\BS(1,m)$
if and only if they are abstractly commensurable, which happens if
and only if $m,n$ are positive integer powers of the same positive
integer.

Theorem A says that every finitely generated group quasi-isometric to
$\BS(1,n)$ can be obtained from $\BS(1,n)$ by first passing to some abstractly
commensurable group and then to some finite extension. We describe this
phenomenon by saying that the group $\BS(1,n)$ is \emph{quasi-isometrically
rigid}. This property is even stronger than what we know for nilpotent groups,
for while Gromov's theorem says that the class of nilpotent groups is a
quasi-isometrically rigid class, outside of a few low-dimensional cases it is
not known whether an individual nilpotent group must always be
quasi-isometrically rigid.

\paragraph*{Comparison with lattices} Recent work on lattices in semisimple Lie
groups has established the quasi-isometric classification of all such lattices. 

In the case of a nonuniform lattice $\Lambda$ in a semisimple Lie
group $G \ne \SL(2,\R)$, quasi-isometric rigidity of $\Lambda$
follows from the deep fact that the quasi-isometry group
$\QI(\Lambda)$ is the commensurator group of $\Lambda$ in $G$, a countable
group (see \cite{Schwartz:RankOne}, \cite{Schwartz:Diophantine},
\cite{FarbSchwartz}, \cite{Eskin:HigherRank}, or
\cite{Farb:Lattices} for a survey).

In contrast, for uniform lattices $\Lambda$ in the isometry group
of $X = \hyp^n$ or $\complex\hyp^n$, the quasi-isometry group
$\QI(\Lambda) \isomorphic \QI(X)$ is $\QC(\bdy X)$, the
(Heisenberg) quasiconformal group of the sphere at infinity $\bdy
X$, an infinite dimensional group. In this situation it is the
whole collection of lattices in $\Isom(X)$ which is
quasi-isometrically rigid, not any individual lattice
\cite{Tukia:quasiconformal}, \cite{CannonCooper},
\cite{Chow:ComplexHyp}.

For $\BS(1,n)$ something interesting happens. It exhibits both
types of contrasting behavior just described: the group $\BS(1,n)$
is quasi-isometrically rigid \emph{and yet} its quasi-isometry
group is infinite-dimensional, as follows. Let $\Q_n$ be the
metric space of ``\nb{$n$}adic rational numbers''. Let $\Bilip(X)$
denote the group of bilipschitz homeomorphisms of the metric
space $X$. In \cite{FarbMosher:BSOne}, Theorem 8.1 we showed:
$$\QI(\BS(1,n)) \homeo \Bilip(\reals) \cross \Bilip(\Q_n) $$ In
proving Theorem A, this formula plays a role similar to that
played by the formula $\QI(X) \homeo \QC(\bdy X)$ when $X=\hyp^n$
or $\complex\hyp^n$.

\vfill\break

\subsection*{Outline of the paper}

As we shall see, the ``mixed behavior'' of $\BS(1,n)$ allows for
some analogies with proof techniques developed in the case of
lattices. Some fundamentally new phenomena occur, however, and
these require new methods.  We point out in particular:
\begin{itemize}
\item The notion of quasisimilarity (a delicate though crucial
variant of quasisymmetric map) and the corresponding notion of
dilatation.
\item The theory of \emph{biconvergence groups}.
\item A method which applies these dynamical properties to proving
quasi-isometric rigidity.
\end{itemize}

\paragraph*{Section \ref{section:background}: Geometry and
boundaries for $\mathbf{BS(1,n)}$.} We review some results of
\cite{FarbMosher:BSOne}. We construct a metric 2-complex $X_n$ on
which $\BS(1,n)$ acts properly discontinuously and cocompactly by
isometries, and we equip this complex with a boundary which is
formed of two disjoint pieces: an \emph{upper boundary} $\bdy^u
X_n \homeo\Q_n$, which is the space of hyperbolic planes in $X_n$,
and a \emph{lower boundary} $\bdy_\ell X_n \homeo \reals$. We also
review the fact that a quasi-isometry of $X_n$ induces bilipschitz
homeomorphisms of the upper and lower boundaries (Proposition
\ref{prop:bilip}), giving the isomorphism $$\QI(\BS(1,n)) \homeo
\Bilip(\reals) \cross \Bilip(\Q_n) $$

\paragraph*{Section \ref{section:quasiaction}: Representation into
the quasi-isometry group.} We recall a principle of Cannon and
Cooper, that if a group $G$ is quasi-isometric to a proper
geodesic metric space $X$, then there is a ``quasi-action'' of $G$
on $X$, and an induced representation $G \to \QI(X)$. Combining
this with the results of \S\ref{section:background} we obtain the
fact that if $G$ is quasi-isometric to $\BS(1,n)$ then there is an
induced action $\rho = \rho_\ell \cross \rho^u \from
G\to\Bilip(\reals) \cross \Bilip(\Q_n)$.

\textsl{Important point:\/} At this point one might try to make an
analogy with \cite{Tukia:quasiconformal} (see also
\cite{CannonCooper}), and attempt to prove that a ``uniformly
quasiconformal'' subgroup of $\Bilip(\reals) \cross \Bilip(\Q_n)$
is conjugate into some kind of ``conformal'' subgroup of
$\Bilip(\reals) \cross \Bilip(\Q_n)$. However, there are serious
difficulties with this approach (see comments in
\S\ref{section:questions}).

\paragraph*{Section \ref{section:quasisim}: Uniform
quasisimilarity actions on $\reals$.} Instead we consider the
projected representation $\rho_l \from G \to \Bilip(\reals)$.  We
show in in Proposition \ref{prop:BilipAction}  that $\rho_l(G)$ is
a group of \emph{quasisimilarities} which is uniform with respect
to a certain dilatation.

Inspired by Hinkkanen's Theorem \cite{Hinkkanen:quasisymmetric},
we prove in Theorem \ref{theorem:qshink} that a uniform group of
quasisimilarities of the real line is bilipschitz conjugate to an
affine group (Theorem \ref{theorem:qshink}).  Applying this to the
representation $\rho_\ell\from G\to\Bilip(\reals)$ we obtain a
bilipschitz conjugate representation $\theta\from G \to
\Aff(\reals)$. It is crucial that this conjugacy is bilipschitz as
opposed to just quasisymmetric.  This is why we cannot use
Hinkkanen's original theorem.

\textsl{Important point:\/} For an arbitrary group quasi-isometric
to a lattice in a semisimple Lie group, finding a representation
into that Lie group usually finishes the proof of quasi-isometric
rigidity. In the present case more work is required. One reason is
that neither $G$ nor even $\BS(1,n)$ is a lattice in $\Aff(\R)$,
although $\BS(1,n)$ is a nondiscrete subgroup of $\Aff(\R)$. In
fact $\Aff(\R)$ is a nonunimodular Lie group and so does not admit
any lattice. Hence we must find another way to prove that $\theta$
has finite kernel, and to analyze the image group $\theta(G)$.

\paragraph*{Section \ref{section:affrep}: Biconvergence groups.} 
We study the action of $G$ on the boundary pair $(\bdy_\ell
X_n,\bdy^u X_n)$. Exploiting analogies with convergence groups,
the dynamical behavior of this pair of actions is encoded in what
we call a \emph{biconvergence group}. Using this boundary
dynamics, we show that the representation $\theta \from G \to
\Aff(\reals)$ is virtually faithful
(Proposition~\ref{proposition:virtfaithful}), and that the group
of affine stretch factors, or \emph{stretch group} of
$\theta(G)\subset \Aff(\R)$ is infinite cyclic
(Proposition~\ref{proposition:cyclic:stretch}). The proof of
Proposition~\ref{proposition:cyclic:stretch} makes vital use of a
\emph{bilipschitz} conjugacy between $\rho_\ell$ and $\theta$.

\paragraph*{Section \ref{section:proof}: Finishing the proof of Theorem A.} We
apply combinatorial group theory and quasi-isometry invariants, together with
the results of \S\ref{section:affrep}, to identify the image group $\Gamma =
\theta(G) \subset \Aff(\R)$. We show that $\Gamma$ is the mapping torus of some
injective, nonsurjective endomorphism $\phi \from A \to A$ where $A$ is either
the infinite cyclic group or the infinite dihedral group---that is, $\Gamma$
has the presentation $\left<A,t \suchthat tat^\inverse = \phi(t), \, \forall a
\in A\right>$. In particular $\Gamma$ contains a subgroup of index $\le 2$
isomorphic to $\BS(1,m)$ for some integer $m \ge 2$. Applying
\cite{FarbMosher:BSOne} it follows that $\BS(1,m)$ is abstractly commensurable
to $\BS(1,n)$, and so the same is true of $\Gamma$. We also prove in Corollary
\ref{corollary:torsionfree} that if $G$ is torsion free then $G$ is isomorphic
to $\BS(1,k)$ for some integer $k$ with $\absvalue{k} \ge 2$.

\paragraph*{Section \ref{section:questions}: Final comments.} We
discuss the possibility of strengthening Theorem A.

\subsection*{Acknowledgements} We thank J. Sturm and M. Feighn for useful
suggestions and discussions, and A. Hinkkanen for helpful correspondence. We
also thank the members of the Workship on Hyperbolic and Automatic Groups at the
University of Montreal where a preliminary version of this work was presented,
especially S. Gersten for helpful criticism, and T. Delzant for suggesting
the theorem of Bieri-Strebel which we use in \S\ref{section:proof}.

\section{Geometry and boundaries for $\mathbf{\BS(1,n)}$}
\label{section:background}

In this section we briefly review the material from \cite{FarbMosher:BSOne}
which we  will need.

\subsection{Quasi-isometries}
A $(K,C)$\emph{quasi-isometry} between metric spaces is a map $f\from
X\to Y$ such that, for some constants $K,C>0$:
\begin{itemize}
\item For all $x_1, x_2 \in X$ we have
$$\frac{1}{K} \, d_X(x_1,x_2) - C \le d_Y(f(x_1),f(x_2)) \le K
d_X(x_1,x_2) + C$$
\item The $C$-neighborhood of $f(X)$ is all of $Y$.
\end{itemize}
A \emph{coarse inverse} of a quasi-isometry $f \from X \to Y$ is a
quasi-isometry $g \from Y \to X$ such that, for some constant $C'>0$, we
have $d(g\circ f(x),x)<C'$ and $d(f\circ g(y),y)<C'$ for all $x \in X$
and $y \in Y$. Every $(K,C)$ quasi-isometry $f \from X \to Y$ has a
coarse inverse $g \from Y \to X$, namely, for each $y \in Y$ define
$g(y)$ to be any point $x \in X$ such that $d(f(x),y) \le C$.
\medskip

A fundamental observation due to Effremovich-Milnor-\Svarc\ states that
if a group $G$ acts properly discontinuously and cocompactly by
isometries on a proper geodesic metric space $X$, then $G$ is finitely
generated, and $X$ is quasi-isometric to $G$ equipped with the word
metric. Two finitely generated groups $G_1,G_2$ are said to be
\emph{(abstractly) commensurable} if there are finite index subgroups
$H_i \subgroup G_i$ such that $H_1,H_2$ are isomorphic to each other. It
is easy to check that abstractly commensurable groups are
quasi-isometric, and that if $G$ is finitely generated and $N$ is a
finite normal subgroup then $G$ is quasi-isometric to $G/N$.

Let $\QIMap(X)$ be the set of all quasi-isometries $f \from X \to X$, equipped
with the binary operation of composition. Given $f,g \in \QIMap(X)$ and $C \ge
0$, we write $f \Hequiv_C g$ if $d(f(x),g(x))<C$ for all $x \in X$. We 
write $f \Hequiv g$ if there exists $C \ge 0$ such that $f \Hequiv_C g$; this is
an equivalence relation on $\QIMap(X)$, known as \emph{Hausdorff
equivalence} or \emph{coarse equivalence}. The set of equivalence classes is
denoted $\QI(X)$. The operation of composition respects Hausdorff equivalence,
in the sense that if $f_1 \Hequiv f_2$ and $g_1 \Hequiv g_2$ then $f_1 \composed
g_1 \Hequiv f_2 \composed g_2$. Composition therefore descends to a well-defined
binary operation on $\QI(X)$.  With respect to this operation, $\QI(X)$ is a
group, whose identity element is the Hausdorff equivalence class of the identity
map on~$X$. Inverses exist in $\QI(X)$ because of the fact that every
quasi-isometry has a coarse inverse.

\subsection{The 2-complex $X_n$}

Throughout this paper we use the upper half plane model of the hyperbolic plane
$$\hyp^2 = \{ (x,y) \in \reals^2 \suchthat y >0 \}
$$
with metric $ds^2 = (dx^2 + dy^2) / y^2$. Also, let $T_n$ denote the unique
homogeneous, directed tree such that each vertex has one incoming directed edge
and $n$ outgoing directed edges, with the geodesic metric that makes each edge
of $T_n$ isometric to the interval $[0,\log(n)]$.

In \cite{FarbMosher:BSOne} we constructed a metric 2-complex $X_n$ on
which $\BS(1,n)$ acts properly discontinuously and cocompactly by
isometries. The 2-complex $X_n$ is homeomorphic to $T_n \cross \reals$.
There is a geodesic metric on $X_n$ with the following properties:
\begin{itemize}
\item For each directed edge $E \subset T_n$, the subset of $X_n$
corresponding to $E \cross \reals$ is isometric to the ``horostrip of
height $\log(n)$'', namely 
$$\{(x,y) \in \hyp^2 \suchthat 1 \le y \le n\}
$$
\item For each bi-infinite directed line $L \subset T_n$, the subset of
$X_n$ corresponding to $L \cross \reals$ is isometric to $\hyp^2$. 
\item For any two bi-infinite directed lines $L \ne L' \subset T_n$,
the subset of $X_n$ corresponding to $(L \intersect L') \cross \reals$
is isometric to a ``horodisc exterior'' $\{(x,y) \in \hyp^2 \suchthat y
\le 1\}$.
\end{itemize}
A cocompact, properly discontinuous, free action of $\BS(1,n)$ on $X_n$
is described in \cite{FarbMosher:BSOne} by exhibiting $X_n$ as the
universal cover of a certain metric 2-complex $C_n$ whose fundamental
group is $\BS(1,n)$.

Various important features of the complex $X_n$, and the action of
$\BS(1,n)$, are summarized in the following proposition. Given a
continuous map between metric spaces $f \from A \to B$, an
\emph{isometric section} is a subset $B' \subset A$ such that $f
\restrict B'$ is an isometry onto $B$. The isometric sections of the
\emph{height function} $h \from \hyp^2 \to \R$, defined by $h(x,y) =
\log(y)$, are precisely the \emph{vertical lines} in $\hyp^2$, lines of
the form $x=\text{(constant)}$. Also define a height function $h \from
T_n \to \R$, by requiring $h(v_0)=0$ for some chosen base vertex $v_0 \in
T_n$, and requiring that $h$ takes each edge of $T_n$ onto some interval
of length $\log(n)$ in $\R$ by an orientation preserving isometry. The
isometric sections of $h\from T_n \to\reals$ are precisely the
bi-infinite directed lines in $T_n$ (called ``coherent lines'' in
\cite{FarbMosher:BSOne}).

\begin{proposition} There exist actions of $\BS(1,n)$ on $X_n$, $\hyp^2$,
$T_n$, and $\reals$, and equivariant maps between these spaces as summarized in
the following commutative diagram:
$$\xymatrix{
& X_n \ar[dl]_p \ar[dd]^h \ar[dr]^q \\
{\hyp^2} \ar[dr]_h & & T_n \ar[dl]^h \\
& {\reals}
}
$$
Moreover:
\begin{itemize}
\item The action on $X_n$ is properly discontinuous, cocompact, and free.
\item The function $h \from X_n \to \R$ is the fiber product of the
height functions $h \from \hyp^2 \to \R$ and $h \from T_n \to
\R$. That is, the map $X_n \to \hyp^2 \cross T_n$ taking $x$ to
$(p(x),q(x))$ is an equivariant homeomorphism onto the subset of $\hyp^2
\cross T_n$ consisting of all pairs $(a,b)$ such that $h(a)=h(b)$.
\item The map $q \from X_n \to T_n$ induces a 1--1 correspondence between
vertical lines in $T_n$ and isometric sections of the map $p\from
X_n\to\hyp^2$.
\item The map $p \from X_n \to \hyp^2$ induces a 1--1 correspondence
between vertical lines in $\hyp^2$ and isometric sections of the map
$q\from X_n \to T_n$.
\end{itemize}
\label{prop:actions}
\end{proposition}

Isometric sections of $p \from X_n \to \hyp^2$ are called \emph{hyperbolic
planes} in $X_n$, and isometric sections of $q \from X_n \to T_n$ are called
\emph{trees in $X_n$}. Since $\BS(1,n)$ clearly acts on the set of vertical
lines in $T_n$, it follows that $\BS(1,n)$ acts on the set of hyperbolic planes
in $X_n$. Similarly, $\BS(1,n)$ acts on the set of trees in $X_n$.

The proof of the above proposition can be gleaned from the information in
\cite{FarbMosher:BSOne}. Here is an alternative proof, which gives an
interesting new construction of $X_n$.

\begin{proof}
Recall the presentation $\BS(1,n) = \left< a,b \suchthat bab^\inverse =
a^n\right>$. Define the \emph{height action} of $\BS(1,n)$ on $\reals$ by
$$
a \cdot t = t, \qquad b \cdot t = t + \log(n)
$$
To describe the actions of $\BS(1,n)$ on $\hyp^2$ and $T_n$, we first define
\emph{affine actions} on $\reals$ and on the $n$-adic rational numbers $\Q_n$,
and then we describe how these induce the desired actions on $\hyp^2$ and $T_n$.

If $R$ is a ring with unit, let $\Aff(R)$ be the group of all matrices of the
form $\textmatrix{r}{s}{0}{1}$, where $r,s \in R$ and $r$ is invertible; the
group law is ordinary matrix multiplication. The group $\Aff(R)$ acts on $R$ by
fractional linear transformations: $\textmatrix{r}{s}{0}{1} \cdot x = rx+s$ for
all $x \in R$. Note that if $R$ is a commutative ring then $\Aff(R)$ is a
solvable group. 

If the integer $n$ is invertible in the ring $R$, there is a
representation $\BS(1,n) \mapsto \Aff(R)$ defined by
\begin{align*}
&a \mapsto \begin{pmatrix}1 & 1 \\ 0 & 1\end{pmatrix} 
&b \mapsto \begin{pmatrix}n & 0 \\ 0 & 1\end{pmatrix} \\
&a \cdot r = r+1 &b \cdot r = nr
\end{align*}
for all $r \in R$. We call this the \emph{affine action} of $\BS(1,n)$
on the ring $R$.

As a special case we obtain an affine action of $\BS(1,n)$ on the real
numbers. The action of $\Aff(\reals)$ on $\reals$ extends to an
isometric action on $\hyp^2$, and by composition we obtain the desired
action of $\BS(1,n)$ on $\hyp^2$.

For another special case, let $\Q_n$ be the ring of $n$-adic rational
numbers. This is the ring of all formal series $\Sum_{i\in\Z} \zeta_i
n^{i}$ with $\zeta_i \in \{0,\ldots,n-1\}$ for all~$i$, and $\zeta_i = 0$
for $i$ sufficiently close to $-\infinity$, with the obvious addition
and multiplication. We write the series $\Sum_{i\in\Z} \zeta_i n^{i}$
more succinctly as $(\zeta_i)_{i \in \Z}$ or just $(\zeta_i)$ when $i
\in \Z$ is understood. The ring of $n$-adic integers $\Z_n$ is the
subring of all $(\zeta_i) \in \Q_n$ with $\zeta_i=0$ for all
$i<0$. Note that the integer $n$ is invertible in $\Q_n$, and so
the affine action of $\BS(1,n)$ on $\Q_n$ is defined.

The metric on $\Q_n$ is the usual \emph{$n$-adic
metric}, where the distance between $(\zeta_i)$ and $(\zeta'_i)$ in
$\Q_n$ is $n^{-k}$, where $k\in \integers$ is the maximum integer such
that $\zeta_i = \zeta'_i$ for all $i\le k$. 

The tree $T_n$ may be identified with the Bruhat-Tits building of
$\Q_n$, and so the action of $\BS(1,n)$ on $\Q_n$ induces an action on
$T_n$. To be explicit, for each truncated series $\eta = (\eta_i)_{i\le
k}$, where $\eta_i \in \{0,\ldots,n-1\}$ for $i \le k$, and $\eta_i = 0$
for $i$ sufficiently close to $-\infinity$, define a set 
$$C_\eta = \{(\zeta_i) \in \Q_n
\suchthat \zeta_i=\eta_i, \, \forall i \le k \}
$$
The set $C_\eta$ is called a \emph{clone of $\Z_n$ in $\Q_n$}. The integer $k$ is
called the \emph{combinatorial height} of the clone $C_\eta$, denoted
$h_c(C_\eta)$. The inclusion lattice on clones defines the directed tree
$T_n$ as follows. The vertices of $T_n$ are the clones of $\Z_n$ in $\Q_n$. There
is a directed edge $C_\eta \to C_{\eta'}$ if and only if the following hold:
\begin{itemize}
\item $C_\eta \supset C_{\eta'}$.
\item For all $C_{\eta''}$, if $C_\eta \supset C_{\eta''} \supset C_{\eta'}$ then
$C_{\eta''} = C_\eta$ or $C_{\eta''}=C_{\eta'}$.
\end{itemize}
Note that $h_c(C_{\eta'}) = h_c(C_\eta) + 1$. Note also that $T_n$ \emph{is} a
tree, because any two clones are either disjoint or one contains the other, and
it is easy to check that each vertex has one incoming and $n$ outgoing edges.

The action of $\Aff(\Q_n)$ on $\Q_n$ takes clones to clones, preserving
inclusion. The affine action of $\BS(1,n)$ on $\Q_n$ therefore induces a
direction preserving action of $\BS(1,n)$ on the tree $T_n$. 

Now we may define $X_n$, and the height function $h \from X_n \to \reals$, by
applying the fiber product construction to the height functions $h\from\hyp^2 \to
\R$ and $h \from T_n \to \R$. Since the latter two height functions are
$\BS(1,n)$ equivariant, we obtain an action of $\BS(1,n)$ on $X_n$ so that $h
\from X_n \to \reals$ is equivariant. 

To see that the action on $X_n$ is properly discontinuous, note that the
general element of $\BS(1,n)$ is a matrix of the form
$\textmatrix{n^i}{k/n^j}{0}{1}$, and if a 1--1 sequence of such matrices
is bounded in $\Aff(\R)$ then it is unbounded in $\Aff(\Q_n)$. The rest
of the proof follow  easily.  
\end{proof}

\begin{Remarks} 
\item This construction of $X_n$ is equivalent to the construction
given in \cite{FarbMosher:BSOne}.
\item The action of $\BS(1,n)$ on $T_n$ is isomorphic to the action on the
Bass-Serre tree of the HNN decomposition $\BS(1,n) \homeo
\integers *_\phi$ where $\phi \from \integers \to \integers$ is given by
$\phi(k)=nk$. 
\item The action of $\BS(1,n)$ on $\hyp^2$ is a ``laminable action'' as
described in \cite{Mosher:indiscrete}. The action on $T_n$ is also laminable:
the decomposition of $X_n$ into trees gives a $T_n$-lamination which resolves
the indiscreteness of the action of $\BS(1,n)$ on $T_n$. 
\item Let $\Z[1/n]$ be the ring of fractions obtained from $\Z$ by
inverting $n$. Note that the affine representation $\BS(1,n) \to
\Aff(\Z[1/n])$ is a group isomorphism. The natural inclusions from
$\Aff(\Z[1/n])$ into $\Aff(\R)$ and $\Aff(\Q_n)$ then give 
the affine representations $\BS(1,n) \to \Aff(\R)$ and
$\BS(1,n)\to\Aff(\Q_n)$.

\end{Remarks}

\subsection{The upper and lower boundaries of $X_n$}
The boundary of the hyperbolic plane in the upper half plane model may be
written as $\bdy\hyp^2 = \reals \union \{+\infinity\}$. Using the
affine action of $\BS(1,n)$ on $\reals$, Proposition
\ref{prop:actions} gives $\BS(1,n)$-equivariant bijections
$$\{\text{trees in } X_n\} \approx \{\text{vertical lines in } \hyp^2\}
\approx \reals
$$
The \emph{lower boundary} $\bdy_\ell X_n$ is defined to be any of these
objects. Note that any two hyperbolic planes $Q, Q' \subset X_n$
intersect in a common horodisc exterior, and so $\bdy Q$ and $\bdy Q'$
share a line at infinity which may be identified with $\bdy_\ell X_n$.
This was how $\bdy_\ell X_n$ was defined in \cite{FarbMosher:BSOne}.

As a dual picture we have the \emph{upper boundary} $\bdy^u X_n$, defined
to be the set of isometrically embedded hyperbolic planes in $X_n$. 
The space of ends of the directed tree $T_n$ can be written as
$\Ends(T_n) = \Q_n \union \{-\infinity\}$ where $\Q_n$ is naturally
identified with the set of positively asymptotic ends and $-\infinity$ is
the unique negatively asymptotic end. By Proposition
\ref{prop:actions} we have $\BS(1,n)$-equivariant bijections
$$\{\text{hyperbolic planes in } X_n\} \approx \{\text{vertical lines in
} T_n \} \approx \Q_n
$$
and so $\bdy^u X_n$ may be identified with any of these.

\begin{Remark}
The set $\overline X_n = X_n \union \bdy_\low X_n \union \bdy^u X_n$ may
be given a topology so that $X_n$ is dense and so that the action of 
$\BS(1,n)$ on $X_n$ extends to a continuous action  on $\overline X_n$.
However $\overline X_n$ is not compact, nor even locally
compact---compare this with the compactifications of
$X_n$ described in \cite{Bestvina:boundaries}.
\end{Remark}

\subsection{Metrics on $\bdy_\low X_n$ and $\bdy^u X_n$}
The metric on $\bdy_\low X_n = \reals$ is the usual metric, and the
metric on $\bdy^u X_n = \Q_n$ is the $n$-adic metric discussed above.

There are a few other ways of visualizing this metric. If $(\zeta_i),
(\zeta'_i) \in \Q_n$ correspond  to hyperbolic planes $Q,Q' \subset
X_n$, and if $S$ is the common horodisc exterior $Q \intersect Q'$, then
$d(Q,Q') = e^{-h(\bdy S)}$. Equivalently, if $L,L'$ are the vertical
lines in $T_n$ corresponding to $\zeta,\zeta'$, and if $v \in T_n$ is
the finite endpoint of the ray $L \intersect L'$, then $d(L,L') =
e^{-h(v)}$. 

Note that $\bdy^u X_n$ is a proper metric space, that is, closed balls are
compact. The Hausdorff dimension of $\Q_n$ equals $1$. The metric $d$ is an
\emph{ultrametric}, also called a \emph{nonarchimedean metric}, in other words a
metric satisfying $d(x,z) \le \Max\{d(x,y),d(y,z)\}$ for any $x,y,z \in \bdy^u
X_n$. Distance in $\Q_n$ takes values in the discrete set $\{n^k \suchthat
k \in \integers\}$. For any $\zeta \in \Q_n$ and any $k \in \integers$, the
closed ball around $\zeta$ of radius $n^{-k}$ is precisely the clone of
combinatorial height $k$ that contains $\zeta$.

\subsection{The groups $\Aff(\Q_n)$ and $\Sim(\Q_n)$}
Define a \emph{similarity} of a metric space $X$ to be a bijection $\phi
\from X \to X$ such that, for some constant $k>0$ we have $d(\phi(x),\phi(y)) =
k \, d(x,y)$ for all $x,y \in X$. The number $k$ is called the \emph{stretch
factor} of $\phi$, denoted $\Stretch(\phi)$; we also say that $\phi$ is a
\emph{$k$-similarity}. The set of all similarities of $X$ forms a group under
composition, denoted $\Sim(X)$. 

The groups $\Sim(\reals)$ and $\Aff(\reals)$ acting on $\reals$ are obviously
identical, but the situation is different in $\Q_n$. 

First note that there is a natural isomorphism between $\Sim(\Q_n)$ and the
group $\Aut(T_n)$ of direction preserving automorphisms of $T_n$. Since the
closed balls in $\Q_n$ are precisely the clones, each similarity of $\Q_n$
takes clones to clones preserving inclusion. Every element of $\Aut(T_n)$ clearly
arises in this manner. 

The group $\Aff(\Q_n)$ acts by similarities on $\Q_n$: given
$M=\textmatrix{\zeta}{\zeta'}{0}{1} \in \Aff(\Q_n)$, if $k$ is the least integer
such that $\zeta_k \ne 0$, then $\Stretch(M) = n^{-k}$. We therefore have a
natural monomorphism $\Aff(\Q_n) \inject \Sim(\Q_n)$, but it is not surjective.
For example, given $x,y \in \Q_n$ such that $x-y$ is invertible in $\Q_n$, the
only affine transformation fixing $x$ and $y$ is the identity, but for any clone
$C$ not containing $x$ or $y$ the subgroup of $\Sim(\Q_n)$ fixing
$x,y$ and preserving $C$ acts transitively on $C$, as can be seen by
constructing appropriate automorphisms of $T_n$.

\subsection{Boundary maps induced by a quasi-isometry}
In \cite{FarbMosher:BSOne} we showed that any quasi-isometry $f \from X_m \to
X_n$ induces maps $f^u \from \bdy^u X_m \to \bdy^u X_n$ and $f_\low \from
\bdy_\low X_m \to \bdy_\low X_n$ characterized as follows. For each hyperbolic
plane $Q \subset X_m$, there is a unique hyperbolic plane $f^u(Q) \subset X_n$
such that 
$$d_\Haus(f(Q),f^u(Q)) < C < \infinity
$$
where $d_\Haus$ denotes Haudorff distance,  and $C$ depends only on the
quasi-isometry constants of $f$. Also, for each hyperbolic plane $Q \subset X_m$,
the closest point projection  from $f(Q)$ to $f^u(Q)$ induces a map from the line
at infinity of $Q$ to the line  at infinity of $f^u(Q)$,  which induces the desired
map $f_\low$. One could also prove  that for any tree $\tau \subset X_m$, there is a
unique tree $f_\ell(\tau) \subset X_n$ such that
$$d_\Haus(f(\tau),f_\ell(\tau)) < C < \infinity
$$

It is proved in \cite{FarbMosher:BSOne} that the 
maps $f^u$ and $f_\low$ are bilipschitz homeomorphisms.  
In fact what we proved is a little stronger, as seen in
the next proposition.

Given a homeomorphism $h \from X \to Y$ between metric spaces, and given $k \ge
1$, we say that $h$ is \emph{$k$-bilipschitz} if
$$\frac{1}{k} \, d(x,y) \le d(f(x),f(y)) \le k \, d(x,y) \quad\text{for
all}\quad x,y \in X
$$
We find it useful to introduce a more precise notion, as follows. Given $0<a<b$,
we say that $h$ is \emph{$[a,b]$-bilipschitz} if 
$$a \cdot d(x,y) \le d(h(x),h(y)) \le b \cdot d(x,y)  \quad\text{for
all}\quad x,y \in X
$$
We use the notation $[a(h),b(h)]$ to denote the \emph{stretch interval}
of $h$, the smallest subinterval of $\reals_+$ such that $h$ is
$[a(h),b(h)]$-bilipschitz.

For example, ``$k$-bilipschitz'' is equivalent to ``$[1/k,k]$-bilipschitz'',
and ``$k$-similarity'' is equivalent to ``$[k,k]$-bilipschitz''.

The next proposition comes from \cite{FarbMosher:BSOne} Proposition 5.3 and
the following remark.

\begin{proposition} For each $K \ge 1$, $C>0$ there exists $L \ge 1$ such that if
$f \from X_m \to X_n$ is a $(K,C)$-quasi-isometry, then there exist $0<a \le b$
such that $b/a < L$, the map $f_\low \from \bdy_\low X_m \to X_n$ is $[a,b]$
bilipschitz, and the map $f^u \from \bdy^u X_m \to \bdy^u X_n$ is
$[1/b,1/a]$-bilipschitz. 
\label{prop:bilip}
\qed\end{proposition}

\paragraph*{Warning.} This does \emph{not} say that the bilipschitz constants of
$f^u,f_\low$ depend only on $K,C$. Even the natural \emph{isometric} action of
$\BS(1,n)$ on $X_n$ does not induce uniformly bilipschitz actions on $\bdy_\low
X_n, \bdy^u X_n$. For example the element $b \in \BS(1,n)$ acts on $\bdy_\low X_n
\homeo \reals$ as a similarity with stretch factor $n$, and so the best
bilipschitz constant of $b^i$ acting on $\bdy_\low X_n$ is $n^i$, although $b^i$
acts isometrically on $X_n$.

\section{Representation into the quasi-isometry group}
\label{section:quasiaction}
In this short section we use a (now standard) principle due essentially
to Cannon-Cooper \cite{CannonCooper}, also indicated by Gromov
\cite{Gromov:ICMAddress}, to prove that if a finitely generated group $G$ is
quasi-isometric to $\BS(1,n)$ then there is an induced representation $G \to
\QI(\BS(1,n)) \homeo
\Bilip(\reals)\cross\Bilip(\Q_n)$. We also study the quality of this
representation.

\subsection{Quasi-isometries and quasi-actions} 
Given a proper, geodesic metric space $X$ and a 
group $G$, a \emph{quasi-action}
of $G$ on $X$ is a map $\psi \from G \to \QIMap(X)$ 
such that, for some constants
$K \ge 1$ and $C \ge 0$, we have:
\begin{itemize}
\item $\psi(g)$ is a $(K,C)$-quasi-isometry, for all $g \in G$.
\item $\psi(\Id) \Hequiv_C \Id_X$.
\item $\psi(g) \composed \psi(h) \Hequiv_C \psi(g h)$, for all $g,h \in
G$.
\end{itemize}
A quasi-action $\psi \from G
\to \QIMap(X)$ evidently induces a homomorphism 
$\Psi \from G \to \QI(X)$, where
$\Psi(g)$ is the Hausdorff equivalence class of $\psi(g)$.

A quasi-action $\psi \from G \to \QIMap(X)$ is \emph{cocompact} if for some $x
\in X$ and $C \ge 0$ the $C$-neighborhood of the set $\{ \psi(g)(x) \suchthat g
\in G\}$ is $X$. Also, $\psi$ is \emph{properly discontinuous} if for each $x,y
\in X$ and $C \ge 0$ the set $\{g \in G \suchthat d(\psi(g)(x),y) \le C)\}$ is
finite.

\begin{proposition}[QI rigidity condition]
\label{proposition:rigidity}
\quad Let $X$ be a proper geodesic metric space, and $\Gamma$ a finitely generated
group. If $f \from \Gamma \to X$ is a quasi-isometry with coarse inverse $\bar f
\from X \to \Gamma$, and if $L_g \from \Gamma \to \Gamma$ denotes left
multiplication by $g \in \Gamma$, then $\psi(g) = f\circ L_g\circ \bar f$ defines
a properly discontinuous, cocompact quasi-action $\psi \from \Gamma \to
\QIMap(X)$.
\end{proposition}

\begin{proof}
It is evident that $\psi$ is a quasi-action. To see that $\psi$ is cocompact,
for any $x \in X$ we have $\{L_g(\bar f(x)) \suchthat g \in G\} = G$, and so
$\{\psi(g)(x) \suchthat g \in G\} = f(G)$, whose $C$ neighborhood equals $X$. 

To see that $\psi$ is properly discontinuous, fix $x,y \in X$ and $C \ge 0$.
Given $g \in G$, if $d(\psi(g)(x),y) \le C$ then $d(\bar
f\composed\psi(g)(x),\bar f(y)) \le KC+C^2$, and since $\bar f
\composed \psi(g) = \bar f \composed f \composed L_g \composed \bar f$ we also
have $d(\bar f \composed \psi(g)(x),L_g(\bar f(y))) \le C$ and so $d(\bar f(y),
g \cdot \bar f(y)) \le KC + C^2 + C$. Since $\Gamma$ is finitely generated,
there are only finitely many $g$ which satisfy this inequality. 
\end{proof}

\begin{Remarks} 
\item If $\psi \from G \to \QIMap(X)$ is a properly discontinuous, cocompact
quasi-action on a proper geodesic metric space $X$, does $G$ have a true action on
$X$ by quasi-isometries? 
\item The passage from a quasi-action $\psi \from G \to \QIMap(X)$ to its
associated representation $\Psi \from G \to \QI(X)$ seems to involve a
loss of information. The multiplicative constant can be recovered,
because if $f \Hequiv g$ and if $f$ is a $(K,C)$-quasi-isometry then $g$
is a $(K,C')$-quasi-isometry for some $C'$. In particular, one can define
the \emph{dilatation} of an element of $\QI(X)$ as the infimum of $K$
such that each representative $f \in \QIMap(X)$ is a
$(K,C)$-quasi-isometry for some $C$.

It is unclear how to recover the additive constant in general, but for
example if $\psi \from G \to \QI(\hyp^3) = \QC(S^2)$ has bounded
dilatation then $\psi(G)$ is induced by a uniformly quasi-isometric
action of $G$ on $\hyp^3$; this follows from \cite{Tukia:quasiconformal}.
\end{Remarks}

\subsection{Consequences for $\mathbf{BS(1,n)}$}

By Proposition \ref{proposition:rigidity}, if 
$G$ is a finitely generated group
quasi-isometric to $\BS(1,n)$ we obtain a quasi-action $\psi \from G \to
\QIMap(X_n)$. Let $\rho \from G \to \QI(X_n)$ be the induced representation.
Applying Proposition \ref{prop:bilip} to $\psi(g)$ 
for each $g \in G$, we obtain:

\begin{proposition} Let $G$ be a finitely generated group quasi-isometric to
$\BS(1,n)$. Then there is an induced representation $\rho \from G \to
\QI(\BS(1,n)) \homeo \Bilip(\reals) \cross \Bilip(\Q_n)$. Furthermore, let
$\rho_\ell \from G \to \Bilip(\reals)$ and $\rho^u \from G \to \Bilip(\Q_n)$ be
the projected actions, and let $[a_\ell(g),b_\ell(g)]$ and $[a^u(g),b^u(g)]$ be
the stretch intervals of $\rho_\ell(g)$, $\rho^u(g)$ respectively. There exists
$L \ge 1$ such that for all $g \in G$ we have
\begin{itemize}
\item $b_\ell(g) / a_\ell(g) \le L$.
\item $b^u(g) / a^u(g) \le L$.
\item $[a_\ell(g),b_\ell(g)] \cdot [a^u(g),b^u(g)] \subset [1/L,L]$.
\end{itemize}
\label{prop:BilipAction}
\end{proposition}

The multiplication of intervals used in this proposition is defined as follows:
if $[a,b]$, $[c,d]$ are subintervals of $\reals_+$, then 
$$[a,b] \cdot [c,d] = \{xy \suchthat x \in [a,b] \text{ and } y \in 
[c,d]\}$$

\section{Uniform quasisimilarity actions on $\reals$}
\label{section:quasisim}

An important part of our proof of Theorem A will be to
understand uniform groups of quasisimilarities of the real line. This entire
section is devoted to such a study, and is based on the following theorem of
Hinkkanen (the terms are defined below). Let $\Aff_+(\R)$ be the index 2
subgroup of $\Aff(\R)$ which preserves orientation of $\R$.

\begin{theorem}[Hinkkanen's Theorem]
A uniform group of orientation preserving quasisymmetric homeomorphisms of $\R$
is quasisymmetrically conjugate to a subgroup of $\Aff_+(\R)$.
\end{theorem}

We cannot use Hinkkanen's Theorem directly because we will need to make use of a
bilipschitz conjugacy, which is not generally produced by a uniformly
quasisymmetric group. In Theorem \ref{theorem:qshink} we recast Hinkkanen's
Theorem in the quasisimilarity setting. This setting has a pedagogical advantage
as well---the technical details of the proof 
are simpler, and we believe that it
is easier to see the geometric ideas underlying Hinkkanen's proof. 

\subsection{A Hinkkanen Theorem for uniform quasisimilarity groups}
Let $X$ be a metric space. 

\begin{Definition}[quasisimilarity] 
A function $f \from X\to X$ is a \emph{$K$-quasisimilarity} if for each
distinct triple $x,y,z\in X$ we have
\begin{align}
\frac{1}{K} \le \frac{d(f(z),f(x)) / d(z,x)}{d(f(y),f(x)) /
d(y,x)} \le K \tag*{$(*)$}
\end{align}
An action of a group $G$ on $X$ is a \emph{uniform quasisimilarity action} if
there exists $K \ge 1$ such that the action of each $g \in G$ is a
$K$-quasisimilarity of $X$. 
\end{Definition}

\begin{Remarks} 
\item Note that a similarity is the same as a 1-quasisimilarity.
\item This property is called ``quasiconformal'' in the appendix to
\cite{FarbMosher:BSOne}, because it implies that if $S$ is a metric sphere in
$X$ then $f(S)$ is nested between two metric spheres $S_1, S_2$ such that
$\radius(S_2) / \radius(S_1)$ is bounded. However, even more is true: the ratios 
$$\frac{\radius(S_2)}{\radius(S)}
\quad\text{and}\quad\frac{\radius(S)}{\radius(S_1)}
$$ 
lie in a fixed interval in $\reals_+$, independent of the original sphere $S$. For
this reason it now seems more appropriate to us to refer to this property as
``quasisimilarity''.
\item An orientation preserving homeomorphism $f \from \reals \to
\reals$ is said to be \emph{$K$-quasisymmetric} if the above inequality
$(*)$ holds whenever $z-y=y-x$. An orientation preserving
$K$-quasisimilarity is therefore $K$-quasisymmetric. The converse is not
true: $K$-quasisimilarities are bilipschitz and hence absolutely
continuous; whereas there exist $K$-quasisymmetric maps which are not
absolutely continuous.
\end{Remarks}

As Cooper notes in \cite{FarbMosher:BSOne}, if $f$ is a $K$-quasisimilarity
then, fixing two points $z,w \in X$, we have
$$ \frac{1}{K^2} \le \frac{d(f(x),f(y)) / d(x,y)}{d(f(z),f(w)) / d(z,w)}
\le K^2 \quad\text{for all } x,y \in X
$$
With $L = d(f(z),f(w)) / d(z,w)$ it follows that $f$ is $K^2 \cdot
\Max\{L,1/L\}$-bilipschitz. However, this conclusion throws some
information away. What one really obtains from this argument is
that 
$$ K^{-2} L \le \frac{d(f(x),f(y))}{d(x,y)} \le K^2 L 
$$ 
and so $f$ is $[a,b]$-bilipschitz with $a=K^{-2}L$ and $b=K^2 L$.
Thus, if $f$ is a $K$-quasisimilarity then $f$ is
$[a,b]$-bilipschitz for some $a,b$ such that $b/a \le K^4$.
Conversely, it is easy to see that if $f \from X \to X$ is
$[a,b]$-bilipschitz with $b/a \le K$, then $f$ is a
$K$-quasisimilarity. Note that there does \emph{not} exist a
constant $C$ depending only on $K$ such that every
$K$-quasisimilarity is $C$-bilipschitz, however a map $f \from \reals \to \reals$ is a $K$-quasisimilarity if and only if there exists a similarity $g \from \reals \to \reals$ such that $f \composed g$ is $\sqrt(K)$-bilipschitz.

Combining these observations with Proposition \ref{prop:BilipAction}, it follows
that any group quasi-isometric to $\BS(1,n)$ has a uniform quasisimilarity
action on $\bdy_\low X_n = \reals$ and on $\bdy^u X_n = \Q_n$. 

Here is our quasisimilarity version of Hinkkanen's Theorem:

\begin{theorem}[Quasisimilarity Hinkkanen's Theorem]
\label{theorem:qshink}
Let $\rho \from G \to \Bilip(\reals)$ be a uniform quasisimilarity action of a
group $G$ on $\reals$. Then there exists $\phi \in \Bilip(\reals)$ such that
$\phi \composed \rho(g) \composed \phi^\inverse \in \Aff(\reals)$ for all $g \in
G$. Moreover the bilipschitz constant for $\phi$ depends only on the
quasisimilarity constant for $\rho$.
\end{theorem}

\begin{Remark}
Hinkkanen's paper considers only groups that preserve orientation
on $\reals$, but his methods indicate an easy way to reduce the
general case to the orientation preserving case. We use these
methods at the end of our proof.
\end{Remark}

\subsection{Proof of Theorem \ref{theorem:qshink}}
The proof of this theorem is an adaptation of Hinkkanen's proof for
quasisymmetric maps. The first part is the following proposition which proves
Theorem \ref{theorem:qshink} in the special case that $G = \integers$ and $G$
preserves orientation. This is analogous to \S 3 of
\cite{Hinkkanen:quasisymmetric}, but the details are quite different in the
current quasisimilarity setting. After this proposition, the remainder of the
proof is an almost verbatim quotation of Hinkkanen's proof, and we will mention
only the highlights. 

\begin{proposition} Let $\{f^n \suchthat n \in \integers\}$ be a nontrivial,
orientation preserving, uniform quasisimilarity action of $\integers$ on
$\reals$. Then one of two things happens:
\begin{itemize}
\item $f$ has no fixed points, the action is uniformly bilipschitz, and $f$ is
bilipschitz conjugate to a translation.
\item $f$ has a single fixed point $p$, the action is not uniformly bilipschitz,
and $f$ is bilipschitz conjugate to multiplication $M_s(x) = sx$, for some $s
\in (0,\infinity) - \{1\}$.
\end{itemize} 
In either case, there is a conjugating map whose bilipschitz constant depends
only on the quasisimilarity constant for $\{f^n\}$.
\label{proposition:cyclic:Hinkkanen}
\end{proposition}

\begin{proof} Let $K'$ be a quasisimilarity constant for each $f^n$, and so the
stretch interval $[a(f^n),b(f^n)]$ satisfies $b(f^n)/a(f^n) \le (K')^4 = K$.

\paragraph*{Case 1: $f$ has no fixed points.} It follows that $f^n$ has no fixed
points for $n \ne 0$. The stretch interval $[a(f^n),b(f^n)]$ must contain $1$,
for otherwise $f^n$ or $f^{-n}$ would be a contraction mapping of $\reals$ which
always has a fixed point. Since $b(f^n)/a(f^n) \le K^4$ it follows that
$[a(f^n),b(f^n)] \subset [1/K^4,K^4]$ and so $f^n$ is $K^4$-bilipschitz for all
$n$.

By replacing $f$ with $f^\inverse$, if necessary, we may assume that $f(x)>x$
for all $x \in \reals$. Let $x_n = f^n(x_0)$ be the orbit of some point $x_0$,
and so 
$$ \ldots x_{-2} < x_{-1} < x_0 < x_1 < x_2 < \ldots
$$

Let $\alpha = x_1 - x_0$. Define $\phi \restrict [x_0,x_1]$ to be the translation $x
\mapsto x - x_0$. Extend $\phi$ to a homeomorphism of $\reals$ as follows: for
each $n \in \integers$ define $\phi \restrict [x_n,x_{n+1}]$ by $\phi(x) =
\phi(f^{-n}(x)) + n\alpha$. Clearly $\phi$ is a homeomorphism of $\reals$ conjugating
$f$ to the map $x \mapsto x+\alpha$. The maps $\phi \restrict [x_0,x_1]$ and $x
\mapsto x + \alpha$ are $1$-bilipschitz. Since $f^{-n}$ is $K$-bilipschitz it follows
that $\phi \restrict [x_n,x_{n+1}]$ is $K$-bilipschitz. 

The proof is completed by applying the following:

\theoremstyle{plain}
\newtheorem*{Rubber}{Rubber band principle}
\begin{Rubber}
Given metric spaces $X,Y$ and a collection of subsets $\A$ covering $X$, suppose
that $\phi \from X \to Y$ is a homeomorphism such that:
\begin{enumerate}
\item For each $A \in \A$ the map $\phi \restrict A \from A \to \phi(A)$ is
$K$-bilipschitz.
\item For any pair $x,y \in X$ there is a sequence $x=x_0,x_1,\ldots,x_m=y$ such
that $d(x,y) = \Sum_{i=1}^m d(x_{i-1},x_i)$, and for each $i=1,\ldots,m$ there
exists $A \in \A$ with $x_{i-1},x_i \in A$.
\item For any pair $\xi,\eta \in Y$, there is a sequence $\xi =
\xi_0,\xi_1,\ldots,\xi_m=\eta$ such that $d(\xi,\eta) = \Sum_{i=1}^m
d(\xi_{i-1},\xi_i)$, and for each $i=1,\ldots,m$ there exists $A \in \A$ with
$\xi_{i-1},\xi_i \in \phi(A)$.
\end{enumerate}
Then $\phi$ is $K$-bilipschitz.
\qed\end{Rubber}

Applying this principle to $\phi$ using $\A = \{ [x_n,x_{n+1}] \suchthat n
\in \integers \}$, it follows that $\phi$ is $K$-bilipshitz.

\paragraph*{Case 2: $f$ has at least one fixed point.} In this case we will prove
that $f$ is bilipschitz conjugate to the multiplication map $M_s(x)=sx$, for
some $s \in (0,\infinity) - \{1\}$. Note that if $s \ne r$ then $M_s$ and $M_r$
are not bilipschitz conjugate, so we have to search carefully to find the correct
value of $s$.

\paragraph*{Step 1: $f$ has a unique fixed point.}
Suppose $f$ has more than one fixed point. Since $f$ is not the identity map,
the fixed point set is a proper, closed subset of $\reals$, and so the non-fixed
point set is an open subset of $\reals$ one of whose components is a finite
interval $(x,y)$. We have $f(x)=x$ and $f(y)=y$, and replacing $f$ by
$f^\inverse$ if necessary we have $f(z)>z$ for all $z \in (x,y)$. Since no point
in $(x,y)$ is fixed it follows that $f^n(z) \to y$ as $n \to \infinity$. Choose
$z = (x+y)/2$, and so $d(x,z)/d(z,y) = 1$. By choosing $n$ sufficiently large,
the ratio 
$$\frac{d(z,y) \cdot d\bigl(f^n(x),f^n(z)\bigr)}{d(x,z) \cdot
d\bigl(f^n(z),f^n(y)\bigr)} =
\frac{d\bigl(x,f^n(z)\bigr)}{d\bigl(f^n(z),y\bigr)}
$$ 
may be made larger than $K'$, violating the fact that $f^n$ is a
$K'$-quasisimilarity.

\bigskip

Let $p$ be the unique fixed point of $f$. 

\paragraph*{Step 2: $p$ is either attracting or repelling.}
Suppose not.   Then either $f(x)>x$ for all $x \ne p$,  or $f(x)<x$ for all $x
\ne p$. We assume the former and derive a contradiction; the latter case is
similar.

If $x<p$ then $f^n(x) \to p$ as $n \to \infinity$, and if $y>p$ then $f^n(y)
\to +\infinity$ as $n \to \infinity$. Taking $x=p-1$ and $y=p+1$, if $n$ is
sufficiently large we may make the ratio
$$ \frac{d(x,p) \cdot d\bigl(f^n(p),f^n(y)\bigr)}{d\bigl(f^n(x),f^n(p)\bigr)
\cdot d(p,y)} = \frac{d\bigl(p,f^n(y)\bigr)}{d\bigl(f^n(x),p\bigr)}
$$
as large as we like, contradicting that $\{f^n\}$ is a uniform quasisimilarity
action.

\bigskip

Replacing $f$ by $f^\inverse$ if necessary, we will assume that $p$ is
repelling, so $f^n(x) \to +\infinity$ if $x>p$, and $f^n(x) \to -\infinity$ if
$x<p$. Under this assumption the stretch factor $s$ will be $>1$.

\paragraph*{Step 3: $f^n$ is \emph{not} uniformly bilipschitz.}
If $f^n$ is $L$-bilipschitz for all $n$, then for any $y>p$ consider the
sequence $f^n(y)$ as $n \to +\infinity$. Then 
$d(p,f^n(y)) = d(f^n(p),f^n(y)) \le
L \cdot d(p,y)$ and so $f^n(y)$ is a bounded 
sequence. Since $y < f(y) < f^2(y) <
\cdots$ it follows that $f^n(y)$ converges to a 
fixed point of $f$ distinct from
$p$, a contradiction.

\paragraph*{Step 4: Finding the expansion constant \boldmath$s$.}
Using the assumptions that $f$ preserves orientation and $p$ is repelling, we
show that there is a unique real number $s>1$ such that $s^n \in
[a(f^n),b(f^n)]$ for all $n$. 

First we note an interesting property of any uniform quasisimilarity action of
a group $G$ on a metric space $X$. Let $[a(g),b(g)]$ be the stretch interval
for the action of $g$ on $X$. Note that the map $g \mapsto [a(g),b(g)]$ satisfies
the following properties:
\begin{itemize}
\item $b(g)/a(g) \le K$ for all $g \in G$.
\item $a(\Id)=b(\Id)=1$.
\item $[a(gh),b(gh)] \subset [a(g)a(h),b(g)b(h)]$ for all $g,h \in G$.
\item $[a(g^\inverse),b(g^\inverse)] = [b(g)^\inverse,a(g)^\inverse]$ for all
$g \in G$.
\end{itemize}

\begin{Definition}[uniform quasihomomorphism]
Given a group $G$, a \emph{uniform quasihomomorphism} from $G$ to $\reals_+$ is
a map which associates to each $g \in G$ an interval $[a(g),b(g)]$ satisfying the
properties listed above, for some constant $K \ge 1$ independent of $g$.
\end{Definition}

Hence to each uniform quasisimilarity action of a group $G$ on a metric space,
there is an associated uniform quasihomomorphism from $G$ to $\reals_+$.

The following lemma, applied to the uniform quasihomomorphism $n \mapsto
[a(f^n),b(f^n)]$, produces the expansion factor $s$:

\begin{lemma} For any uniform quasihomomorphism $n \mapsto [a_n,b_n]$ from
$\integers$ to $\reals_+$, there exists a unique $s > 0$ such that $s^n \in
[a_n,b_n]$.
\end{lemma}

\begin{proof} Multiplying $\integers$ by $-1$ if necessary, and assuming that
$[a_n,b_n]$ is not uniformly bounded, we may assume that if $n$ is sufficiently
large then $[a_n,b_n] \subset (1,\infinity)$. In this case the $s$ that we find
will be larger than $1$. Let $K$ be a quasihomomorphism constant 
for $n \mapsto [a_n,b_n]$.

For all $n \ge 1$ define the following subinterval of $\reals_+$:
$$I_n = \bigl[\sqrt[n]{\vphantom{b_n} a_n},\sqrt[n]{b_n}\bigr].
$$
Given $m,n \ge 1$ we have $I_{mn} \subset I_n$, because
$$
[a_{mn},b_{mn}] = [a_{\underbrace{n + \cdots + n}_{m \text{ times}}},
b_{\underbrace{n + \cdots + n}_{m \text{ times}}}] 
\subset [a_n^m,b_n^m]
$$
In particular, $I_n \subset I_1 = [a_1,b_1]$ for all $n \ge 1$. Also,
the ratio of the upper and lower endpoints of $I_n$ is $\sqrt[n]{b_n/a_n} \le
\sqrt[n]{K}$, which approaches $1$ as $n \to +\infinity$. Since $I_n$ is a
subinterval of the fixed interval $[a_1,b_1]$ it follows that $\Length(I_n)$
approaches $0$ as $n \to +\infinity$. 

We have shown that the intervals $I_n$ are nested with respect to the
divisor lattice on the natural numbers. Since $\Length(I_n) \to 0$ it
follows that $\Intersection_n I_n$ is a singleton $\{s\}$. This proves the
existence of a unique $s > 1$ such that $s^n \in [a_n,b_n]$ for all $n \in
\integers_+$. 
 
Since $[a_{-n},b_{-n}] = [b_n^\inverse,a_n^\inverse]$, it follows that $s^n \in
[a_n,b_n]$ for all $n \in \integers$.
\end{proof}

\paragraph*{Step 5: $f$ is bilipschitz conjugate to $M_s$.}
Conjugating by a translation if necessary, we may assume that the fixed point of
$f$ is $p=0$. We'll define the conjugacy $\phi$ on $[0,\infinity)$, and prove
it is bilipschitz there. The extension to $(-\infinity,0]$ is defined
similarly, and the rubber band principle proves that $\phi$ is bilipschitz on
$\reals$.

Let $x_0 = 1$, $x_n = f^n(x_0)$. Define $\phi \restrict [x_0,x_1]$ to be the
unique orientation preserving affine homeomorphism from $[x_0,x_1]$ to $[1,s]$.
This map has constant derivative $(s-1)/(x_1-x_0) = (s-1)/(x_1-1)$, and we
must get a bound on this derivative depending only on $K$.

\begin{lemma}
With the notation above,
$$ \frac{s-1}{x_1-1} \in \bigl[ \frac{1}{K^2}, K^2 \bigr]
$$
\label{lemma:phi:bilip}
\end{lemma}

Accepting this lemma for the moment, let us define the conjugacy $\phi$ on all of
$[0,\infinity)$ and prove that it is bilipschitz there, with a constant depending
only on $K$.

Define $\phi(0)=0$, and define $\phi \restrict [x_n,x_{n+1}]$ by
$$\phi(x) = \phi(f^{-n}(x)) \cdot s^n
$$
Obviously $\phi$ is a homeomorphism of $[0,+\infinity)$ conjugating $f$ to
$M_s$. 

We prove that $\phi \from [0,+\infinity) \to [0,+\infinity)$ is bilipschitz by
applying the rubber band principle to the family of sets
$$S_n = \{0, x_n\}, \quad I_n = [x_n,x_{n+1}], \quad n \in \integers
$$
Clearly this family of sets covers $[0,+\infinity)$, and it satisfies
properties (2) and (3) in the rubber band principle. 

Note that 
\begin{align*}
d(0,x_n) &= d(f^n(0),f^n(1)) \\
         &\in [a(f^n),b(f^n)] \cdot d(0,1) \\
         &= [a(f^n),b(f^n)]
\end{align*} 
According to Step 4 this interval is a subset of $[s^n/K, s^n \cdot K]$. Since 
$$d(\phi(0),\phi(x_n)) = d(0,s^n) = s^n
$$
it follows that $\phi$ is $K$-bilipschitz on the set $S_n$. 

To prove that $\phi$ is bilipschitz on $I_n$, we have the following
decomposition of $\phi \restrict I_n$:
$$
\xymatrix{
I_n \ar[rr]_-{f^{-n}} \ar@/^2pc/[rrrrr]^{\phi} & &
[x_0,x_1] \ar[r]_{\phi} &
[1,s] \ar[rr]_-{y \to ys^n} & &
[s^n,s^{n+1}]
}
$$
The first map is $[a(f^{-n}),b(f^{-n})]$-bilipschitz, and this interval is
contained in $[s^{-n}/K, s^{-n} K]$. The second map is $[1/K^2,K^2]$ bilipschitz
by Lemma \ref{lemma:phi:bilip}. The third map is $[s^n,s^n]$-bilipschitz. The
composition is therefore $[1/K^3,K^3]$-bilipschitz. 
\end{proof}

\begin{proof}[Proof of Lemma \ref{lemma:phi:bilip}]
The idea of the proof is that $x_1 / s = d(f(1),f(0))/s \in [a(f),b(f)] \cdot
d(0,1) / s \subset [s/K,K] \cdot 1/s = [1/K,K]$, and so $x_1/s$ is bounded by a
constant depending only on $K$. If $s$ were large it would follow that
$(x_1-1) / (s-1)$ is bounded. However, since $s$ may not be large, we must work
a little harder to make use of the fact that $x_n/s^n$ is also
bounded, by a constant depending only on $K$.

It follows from Step 4 that for all $n \in \integers$ the map $f^n$ is
$[s^n/K,s^n K]$-bilipschitz. In particular,
$$d(x_{n+1},x_n) = d(f^n(x_1),f^n(x_0)) \in [1/K,K] \cdot s^n \cdot d(x_1,x_0)
$$
and so
\begin{align*}
\frac{d(x_{n+1},x_n)}{s^{n+1} - s^n} &\in 
[1/K,K] \cdot \frac{s^n}{s^{n+1}-s^n} \cdot d(x_1,x_0) \\
&= [1/K,K] \cdot \frac{d(x_1,x_0)}{s-1}
\end{align*}
We also have
\begin{align*}
 d(x_n,x_0) & =  d(x_n,x_{n-1}) + \cdots + d(x_1,x_0)                         \\
   & \in [1/K,K] \cdot \frac{d(x_1,x_0)}{s-1} \cdot (s^n - s^{n-1}) + \cdots \\ 
   &\qquad\qquad\qquad + [1/K,K] \cdot \frac{d(x_1,x_0)}{s-1} \cdot (s^1 -
s^0) \\
   & =  [1/K,K] \cdot \frac{d(x_1,x_0)}{s-1} \cdot (s^n - 1)
\end{align*}
and so
$$ \frac{d(x_n,x_0)}{s^n - 1} \in [1/K,K] \cdot \frac{d(x_1,x_0)}{s-1}
$$
On the other hand, we have
$$ x_n = x_n - 0 = f^n(x_0) - f^n(0) \in [1/K,K] \cdot s^n \cdot (x_0 - 0)
$$
and, since $x_0=1$ it follows that
$$\frac{x_n}{s^n} \in [1/K,K]
$$

As $n \to \infinity$, the ratio 
$$ \frac{d(x_n,x_0)}{s^n-1} \biggm/ \frac{x_n}{s^n} = \frac{x_n-1}{s^n-1} \biggm/
\frac{x_n}{s^n}
$$
approaches $1$. Hence, for any $\epsilon>0$, if $n$ is sufficiently large we
have
$$ \frac{d(x_n,x_0)}{s^n - 1} \in [1/(K+\epsilon),K+\epsilon]
$$

We have shown that the two intervals 
$$ [1/K,K] \cdot \frac{d(x_1,x_0)}{s-1} \qquad\text{and}\qquad
[1/(K+\epsilon),K+\epsilon]
$$
have nonempty intersection, both containing $d(x_n,x_0) / (s^n-1)$, and it
follows that
$$ \frac{1}{K+\epsilon} \le K \frac{d(x_1,x_0)}{s-1} \qquad\text{and}\qquad
\frac{1}{K} \cdot \frac{d(x_1,x_0)}{s-1} \le K+\epsilon
$$
and so
$$ \frac{1}{K(K+\epsilon)} \le \frac{d(x_1,x_0)}{s-1} \le K(K+\epsilon)
$$
Since $\epsilon>0$ is arbitrary, the lemma follows.
\end{proof}

Now we prove Theorem \ref{theorem:qshink} in the special case that $G$
preserves orientation of $\reals$, in which case we conjugate $G$ into
$\Aff_+(\reals)$. The proof proceeds along the lines laid out on p. 332 of
\cite{Hinkkanen:quasisymmetric}, with a few comments needed to translate the quasisymmetric setting to
the quasisimilarity setting.

Let $\Homeo_+(\reals)$ be the topological group of orientation preserving
homeomorphisms of $\reals$ in the topology of uniform convergence on compact
sets. Given a sequence $g_i \in \Homeo_+(\reals)$ of $K$-quasisimilarities of
$\reals$, the following two statements are equivalent:
\begin{itemize}
\item The sequence $g_i$ has a convergent subsequence in $\Homeo_+(\reals)$. 
\item There exist $x_1,x_2 \in \reals$ such that $g_i(x_1)$ is bounded and
$|g_i(x_2) - g_i(x_1)|$ is bounded away from $0$ and $\infinity$.
\end{itemize}
In this case, the limiting homeomorphism is also a $K$-quasisimilarity. To
prove that the second statement above implies the first, one uses the
Ascoli-Arzela Theorem, observing that the second statement is equivalent to
saying that $\{g_i\}$ is uniformly bilipschitz and $\{g_i(x_1)\}$ is
bounded for some $x_1$. 

If $G \subset \Homeo_+(\reals)$ is a group of $K$-quasisimilarities of $\reals$,
it follows that the closure $\overline G$ in $\Homeo_+(X)$ is also a group of
$K$-quasisimilarities. If $\overline G$ is bilipschitz conjugate to a subgroup
of $\Aff_+(\reals)$ then so is $G$. We may therefore assume for the rest of the
proof that $G$ is closed in $\Homeo_+(\reals)$.

Lemma 6 of \cite{Hinkkanen:quasisymmetric} says that if $g,h \in G - \{\Id\}$ have distinct fixed
points then $H = g \composed h \composed g^\inverse \composed h^\inverse$ has no
fixed points. Lemma 7 of \cite{Hinkkanen:quasisymmetric} says that if $g,h \in G$ each have no fixed
points and $g \ne h$ then $g(x) \ne h(x)$ for all $x \in \reals$. We need not
reprove these lemmas, because in each case the hypothesis requires only that the
group generated by $g,h$ be uniformly quasisymmetric, which is true in the
present situation where $g,h$ generate a uniform quasisimilarity group. If one
desires, the proofs of Lemmas 6 and 7 may by improved by using the
quasisimilarity conjugacies of Proposition \ref{proposition:cyclic:Hinkkanen},
instead of the quasisymmetric conjugacies of Hinkkanen's paper. 

Define $G_T$ to be the set of all elements of $G$ without fixed points. By the
previous paragraph, $G_T$ is a subgroup of $G$. By Proposition
\ref{proposition:cyclic:Hinkkanen}, $G_T$ consists of those elements which are
bilipschitz conjugate to translations. 

Lemma 8 of \cite{Hinkkanen:quasisymmetric} proves that if $G$ is not cyclic and
$G \ne G_T \ne \{Id\}$ then $G_T$ is not cyclic. Lemma 8 also proves, from the
fact that $G$ is closed in $\Homeo(\reals)$, that $G_T$ is also closed. Lemma 9
of \cite{Hinkkanen:quasisymmetric} proves that if $G_T = \Transl(\reals)$, the
group of all translations of $\reals$, then $G \subset \Aff_+(\reals)$.  

Lemma 10 of \cite{Hinkkanen:quasisymmetric} says, assuming $G$ is closed in
$\Homeo(\reals)$, that the theorem is true in two special cases: $G_T=G$; and
$G_T = \{ \Id \}$. Hinkkanen's proof gives a quasisymmetric conjugacy, but if
one follows the proof through verbatim one sees in the present setting that the
conjugating map is bilipschitz. The point is this. Suppose $G=G_T$. One
constructs a strictly increasing sequence of cyclic subgroups $G_1 \supset G_2
\supset \cdots$ whose union is dense in $G_T$, and then one conjugates each
$G_i$ to a translation group by applying the cyclic version of the theorem,
using a conjugating map $f_i$. Then one shows that there is a limiting map $f =
\lim f_i$ which conjugates all of $G_T$ to a translation group.  In our present
case, the maps $f_i$ are provided by Proposition
\ref{proposition:cyclic:Hinkkanen}, and $f_i$ has a bilipschitz constant
depending only on the quasisimilarity constant of the subgroups $G_i$. Since
each $G_i$ is a subgroup of $G_T$, the quasisimilarity constants of $G_i$ are
uniformly bounded, and so the $f_i$ are uniformly bilipschitz. The limiting map
$f$ is therefore bilipschitz.

The only remaining case is when $G \ne G_T \ne \{\Id\}$, and in this case one
applies Lemmas 8, 9, and 10 as on page 332 of \cite{Hinkkanen:quasisymmetric} to finish the proof in
the case that $G$ preserves orientation.

\bigskip

To complete the proof we consider the case that $G$ does not preserve
orientation. Let $G_0$ be the orientation preserving subgroup of $G$, a normal
subgroup of index $2$. Applying the orientation preserving case of the theorem,
we may do a bilipschitz conjugacy of $G_0$ into $\Aff_+(\reals)$. Applying this
conjugacy to all of $G$, we may assume that $G_0 \subset \Aff_+(\reals)$.
Choose a single element $g \in G - G_0$. 

The most interesting subcase (and the only case we really care about) is when
$(G_0)_T$ is a dense subgroup of $\Transl(\reals)$. We follow the proof of Lemma
9 of \cite{Hinkkanen:quasisymmetric}. There is a dense additive subgroup $S
\subset \reals$, and an isomorphism $\phi \from S \to S$, such that $g(x)+s =
g(x+\phi(s))$ for all $x \in \reals$ and $s \in S$. Since $g$ reverses
orientation on $\reals$, it follows that $\phi$ reverses order on $S$, and so
$\phi$ is continuous on $S$. The map $\phi$ therefore extends to a continuous
isomorphism of the additive group $\reals$. It follows that $\phi(s) = cs$ for
some $c \ne 0$. We therefore have $g(x)+s=g(x+cs)$ for all $x,s \in \reals$.
Taking $x=0$ and $r=cs$ we have $g(r)=c^\inverse r + g(0)$, and so $g$ is
affine. 

If $(G_0)_T$ is not a dense subgroup of $\Transl(\reals)$, then as we have seen
there are two further subcases. In one subcase $G_0$ is a cyclic group of
translations, and in the other case $(G_0)_T$ is trivial and $G_0$ is a group of
dilations with a global fixed point. In particular $G_0$ is abelian and so
$G$ is virtually abelian (these subcases will therefore not arise in
applications to groups quasi-isometric to $\BS(1,n)$, because of Proposition
\ref{proposition:virtfaithful}). The interested reader can show in either of
these cases that one can do a further conjugation keeping $G_0$ in
$\Aff_+(\reals)$ and taking $g$ to an orientation reversing affine
transformation.

\section{A virtually faithful affine action on $\reals$}
\label{section:affrep}

Recall the similarity groups $\Sim(\reals)=\Aff(\reals)$ and $\Sim(\Q_n)$. There
is an obvious inclusion $\Isom(X_n) \inject \Aff(\reals) \cross \Sim(\Q_n)$,
whose image is the subgroup of pairs $(f,g) \in \Aff(\reals) \cross \Sim(\Q_n)$
such that $\Stretch(f) \cdot \Stretch(g) = 1$. We have natural inclusions
\begin{align*}
\BS(1,n) &\inject \Isom(X_n) \\
         &\inject \Aff(\R) \cross \Sim(\Q_n) \\
         &\inject \Bilip(\R) \cross \Bilip(\Q_n)\\
         &\homeo \QI(\BS(1,n))
\end{align*}
the final isomorphism being Theorem 7.1 in \cite{FarbMosher:BSOne}.

There is a split short exact sequence
$$1 \to \Isom(\reals) \to \Aff(\reals) \to \Stretch(\reals) \to 1
$$
where $\Isom(\reals)$ is the group of isometries of $\reals$, and
$\Stretch(\reals)$ is the group of orientation preserving stretch homeomorphisms
of $\reals$, otherwise known as $\SL(1,\reals)$. Note that there are canonical
isomorphisms 
$$\Isom(\reals) \homeo \reals \semidirect (\Z / 2\Z) \quad\text{and}\quad
\Stretch(\reals) \homeo \reals_+ \mapright{\log} \reals
$$
where $\Z / 2\Z$ acts on $\reals$ by a reflection.

Given a subgroup $\Gamma \subgroup \Aff(\reals)$ we therefore obtain a
commutative diagram
$$
\xymatrix{
1 \ar[r] & A(\Gamma) \ar@{^{(}->}[r] \ar@{^{(}->}[d] & \Gamma \ar[r]
\ar@{^{(}->}[d] & \Stretch(\Gamma) \ar[r] \ar@{^{(}->}[d] & 1 \\
1 \ar[r] & \Isom(\R) \ar@{^{(}->}[r] & \Aff(\R) \ar[r] & \Stretch(\R) \ar[r] & 1
}
$$

%

The purpose of this section is to prove:

\begin{theorem} 
\label{theorem:homomorphism}
Let $G$ be a finitely generated group quasi-isometric to $\BS(1,n)$. There exists
a representation $\theta \from G \to \Aff(\reals)$ such that $\ker(\theta)$ is
finite. Moreover, setting $\Gamma = \theta(G)$, the group $\Stretch(\Gamma)$ is
infinite cyclic.
\end{theorem}

To begin the proof, by Proposition \ref{prop:BilipAction} we obtain a properly
discontinuous, cocompact quasi-action $\psi \from G \to \QIMap(X_n)$, inducing a
representation $\rho \from G \to \QI(X_n) \homeo \Bilip(\reals) \cross
\Bilip(\Q_n)$, with projected representations $\rho_\ell \from G \to
\Bilip(\reals)$ and $\rho^u \from G \to \Bilip(\Q_n)$. The representations
$\rho_\ell$ and $\rho^u$ are uniform quasi-similarity actions.

Applying Theorem \ref{theorem:qshink} to the representation
$\rho_\ell$, there exists a bilipschitz homeomorphism $\phi \from \reals
\to \reals$ such that $\theta(g) = \phi \composed \rho_\ell(g) \composed
\phi^\inverse \in \Aff(\reals)$ for all $g \in G$, and so we obtain a
representation 
$$\theta \from G\to \Aff(\R)
$$
We now need to use further properties of the situation to show that the
representation $\theta$ has finite kernel and that the group $\theta(G) \subset
\Aff(\reals)$ has infinite cyclic stretch group.

\subsection{Biconvergence groups}

Recall that if $\Gamma$ is a word hyperbolic group then the action of
$\Gamma$ on its boundary $\bdy\Gamma$ is a \emph{uniform convergence group} 
action, which means that $\Gamma$ acts properly discontinuously and
cocompactly on the triple space 
$$T(\bdy\Gamma,\bdy\Gamma,\bdy\Gamma) = \{(x,y,z) \in \bdy\Gamma \suchthat x \ne
y, y \ne z, z \ne x \}
$$
Furthermore if $G$ is quasi-isometric to $\Gamma$ then the induced action of
$G$ on $\bdy\Gamma$ is also a uniform convergence group action.

Despite the fact that the groups $\BS(1,n)$ are far from being word hyperbolic,
the upper and lower boundaries interact in a way reminiscent of
convergence group actions, motivating the following definition.

\begin{Definition}[Biconvergence group]
Suppose that $X,Y$ are topological spaces. Define two triple spaces
$$T(X,X,Y) = \{(x,y,\zeta) \in X \cross X \cross Y \suchthat
x \ne y \}
$$
and
$$T(X,Y,Y) = \{(x,\eta,\zeta) \in X \cross Y \cross Y \suchthat \eta \ne \zeta\}
$$
A \emph{biconvergence action} of a group $G$ on the pair 
$(X,Y)$ consists of an
action of $G$ on $X$ and an action on $Y$, such that the 
induced diagonal actions
of $G$ on $T(X,X,Y)$ and on $T(X,Y,Y)$ are properly discontinuous. 
If furthermore
the actions of $G$ on $T(X,X,Y)$ and on $T(X,Y,Y)$ are cocompact, 
then we call this a
\emph{uniform biconvergence action} of $G$ on $(X,Y)$.
\end{Definition}

\begin{Example} $\BS(1,n)$ is a uniform biconvergence group on
$(\R,\Q_n)$. Note that $\BS(1,n)$ does \emph{not} act properly discontinuously
on the triple space $T(\R,\R,\R) = \{(x,y,z) \in \R \cross \R \cross \R \suchthat
x\ne y, y\ne z, z \ne x\}$, nor on $T(\Q_n,\Q_n,\Q_n)$. For example the action
on $\R$ contains translations of $\R$ by any rational number of the form $i/n^j$.
\end{Example}

\begin{Example} Let $\Solv$ be the unique 3-dimensional, solvable,
non-nilpotent, connected Lie group which has a cocompact lattice. $\Solv$ can be
identified with $\reals^3$ in such a way that the left invariant metric is 
$$ds^2 = e^{-2t} dx^2 + e^{2t} dy^2 + dt^2
$$
There is an ``upper boundary'' $\bdy^u\Solv$ consisting of the set of all
``right side up'' hyperbolic planes of the form $y=$(constant), and we have an
identification $\bdy^u\Solv = \reals$. There is also a ``lower boundary''
$\bdy_\low\Solv$ consisting of all ``upside down'' hyperbolic planes
$x=$(constant), also identified with $\reals$. If $\Gamma$ is the fundamental
group of a compact 3-manifold $M$ fibering over $S^1$ with torus fiber $T^2$, so
that the monodromy map $T^2 \to T^2$ is an Anosov homeomorphism, then $M$ has a
Riemannian metric locally modelled on $\Solv$, and so $\Gamma$ may be identified
with a cocompact, discrete subgroup of $\Solv$. Under these conditions, one can
show that the induced pair of actions of $\Gamma$ on
$(\bdy_\low\Solv,\bdy^u\Solv) \homeo (\R,\R)$ is a biconvergence action.
\end{Example}

\begin{proposition}[$G$ is a biconvergence group] 
\label{proposition:biconvergence}
Let $G$ be a finitely generated group which is quasi-isometric to $\BS(1,n)$.
The induced action of $G$ on $(\bdy_\low X_n,\bdy^u X_n) \homeo (\reals,\Q_n)$ is
a uniform biconvergence action.
\end{proposition}

\begin{proof}
We use the symbols $\bdy_\low$, $\bdy^u$ as shorthand for $\bdy_\low X_n$,
$\bdy^u X_n$. First we give a detailed proof that the action of $G$ on
$T(\bdy_\low,\bdy_\low,\bdy^u)$ is properly discontinuous and
cocompact.  We then indicate the changes needed for 
$T(\bdy_\low,\bdy^u,\bdy^u)$. 

By Proposition \ref{proposition:rigidity} we have a $(K,C)$-quasi-action $\psi
\from G \to \QIMap(X_n)$ which induces the given action of $G$ on
$T(\bdy_\low,\bdy_\low,\bdy^u)$. 

There is a map $\pi \from T(\bdy_\low,\bdy_\low,\bdy^u) \to X_n$ defined as
follows (see Figure~\ref{FigureBarycenter}). Given $(x,y,\zeta) \in
T(\bdy_\low,\bdy_\low,\bdy^u)$, let $H \subset X_n$ be the hyperbolic plane
corresponding to $\zeta \in \bdy^u$. The boundary $\bdy H$ is identified with
$\bdy_\ell$ plus a point denoted $+\infinity$. Consider the ideal triangle
$\Delta = \Delta(x,y,+\infinity) \subset H$. Define $\pi(x,y,\zeta)$ to be the
barycenter of $\Delta$, that is, the intersection of the perpendiculars from
each vertex of $\Delta$ to the opposite side. 

\begin{figure}
\centeredepsfbox{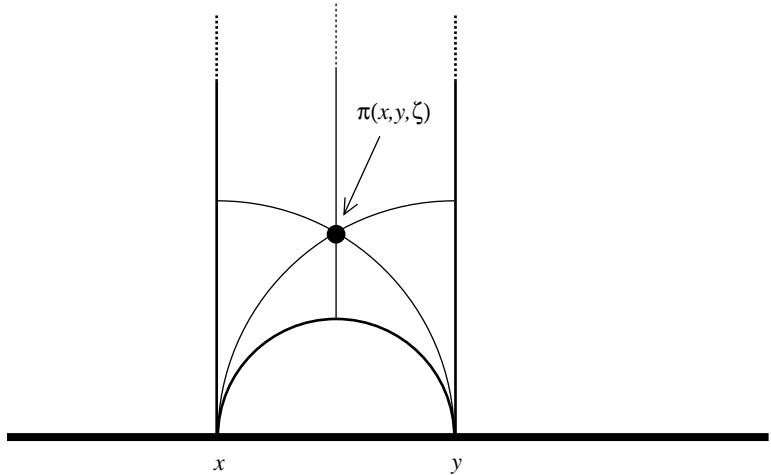}
\caption{The point $\pi(x,y,\zeta)$ is the barycenter of the triangle
$\Delta(x,y,+\infinity)$ in the hyperbolic plane $H \subset X_n$ correponding to
$\zeta$.}
\label{FigureBarycenter}
\end{figure}

The map $\pi$ is obviously continuous, proper, and equivariant with respect to
each isometry $\phi$ of $X_n$, that is
$$\pi(\phi(x,y,\zeta)) = \phi(\pi(x,y,\zeta))
$$

We claim that $\pi$ is almost equivariant with respect to a quasi-isometry. That
is, there exists $C_1$ depending on $K,C$ such that if $f$ is a
$(K,C)$-quasi-isometry of $X_n$, with induced boundary maps $f_\low \from
\R \to \R$ and $f^u \from \Q_n \to \Q_n$, then for all $x,y \in \R, \zeta \in
\Q_n$ we have
$$d\bigl(\pi(f_\low(x),f_\low(y),f^u(\zeta)),f(\pi(x,y,\zeta))\bigr) \le C_1
$$
To prove this claim, note first that $d_\Haus(f(\zeta),f^u(\zeta)) \le C_2$ for
some constant $C_2$ depending only on $K,C$. Composing $f \restrict \zeta$ with
the closest point projection from $f(\zeta)$ to $f^u(\zeta)$, we obtain a
quasi-isometry from the hyperbolic plane $\zeta$ to the hyperbolic plane
$f^u(\zeta)$. This quasi-isometry takes $x,y$, regarded as elements of the line
at infinity of $\zeta$, to $f_\ell(x), f_\ell(y)$, regarded as elements of the
line at infinity of $f^u(\zeta)$, and it takes the point at infinity of $\zeta$
to the point at infinity of $f^u(\zeta)$. The claim now follows from the fact
that the barycenter map of $\hyp^2$ is almost equivariant with respect to
quasi-isometries of $\hyp^2$.

Suppose that the action of $G$ on $T(\bdy_\low,\bdy_\low,\bdy^u)$ is not
properly discontinuous. Then there is a sequence of distinct elements
$g_i \in G$ and a point $(x_0,y_0,\zeta_0) \in T(\bdy_\low,\bdy_\low,\bdy^u)$
such that, setting $(x_i,y_i,\zeta_i) = g_i(x_0,y_0,\zeta_0)$, the sequence
$(x_i,y_i,\zeta_i)$ converges to some $(x,y,\zeta) \in
T(\bdy_\low,\bdy_\low,\bdy^u)$. After chopping off an initial subsequence, it
follows that $\psi_{g_i} (\pi(x_0,y_0,\zeta_0))$ stays in a bounded neighborhood
of $\pi(x,y,\zeta)$ in $X_n$. However, this contradicts the fact that the
quasi-action $\psi$ is properly discontinuous.

Next we prove that the action of $G$ on $T(\bdy_\low,\bdy_\low,\bdy^u)$ is
cocompact. Note that the map $\pi \from T(\bdy_\low,\bdy_\low,\bdy^u) \to X_n$ is
continuous and proper. By cocompactness of the quasi-action $\psi \from G \to
\QIMap(X_n)$, there exists $x_0 \in X_n$ and $C_2>0$ such that for each $y \in
X_n$ there exists $g \in G$ with $d(y,\psi_{g}(x_0)) \le C_2$. 

We claim that the compact set $\pi^\inverse\bigl(\overline B(x_0,K C_2 + 2C +
C_1)\bigr)$ is a fundamental domain for the action of $G$ on
$T(\bdy_\low,\bdy_\low,\bdy^u)$.

To prove this claim, consider any $t \in T(\bdy_\low,\bdy_\low,\bdy^u)$, and
choose $g$ so that $d(\pi(t),\psi_{g}(x_0)) \le C_2$. Applying
$\psi_{g^\inverse}$ we obtain
\begin{align*}
d\bigl( \psi_{g^\inverse}(\pi(t)), \psi_{g^\inverse}
\psi_{g^{\vphantom{\inverse}}} (x_0) \bigr)  &\le K C_2 + C \\
d\bigl( \psi_{g^\inverse}(\pi(t)), x_0) \bigr) &\le K C_2 + C + C
\end{align*}
Also, 
\begin{gather*}
d\bigl(\pi(g^\inverse \cdot t), \psi_{g^\inverse}(\pi(t))\bigr) \le C_1 \\
d(\pi(g^\inverse \cdot t),x_0) \le K C_2 + 2C + C_1 \\
g^\inverse \cdot t \in \pi^\inverse\bigl(\overline B(x_0,K C_2 + 2C +
C_1)\bigr)
\end{gather*}
completing the proof of cocompactness.

To prove that the action of $G$ on $T(\bdy_\low,\bdy^u,\bdy^u)$ is properly
discontinuous and cocompact, it suffices to describe a continuous, proper,
almost $G$-equivariant map $\kappa \from T(\bdy_\low,\bdy^u,\bdy^u) \to X_n$
(see Figure \ref{FigureTreeBarycenter}), and mimic the above proof. Given
$(x,\eta,\zeta) \in T(\bdy_\low,\bdy^u,\bdy^u)$, the point $x$ corresponds to
some tree $\tau \subset X_n$, the image of an isometric section of the projection
map $q \from X_n \to T_n$. The set $\Ends(\tau)$ is identified with $\bdy^u$
plus a point denoted $-\infinity$. The three points
$-\infinity,\eta,\zeta \in \Ends(\tau)$ determine an ideal triangle
$\Delta(-\infinity,\eta,\zeta) \subset \tau$ whose three sides intersect in a
unique point which we take to be $\kappa(x,\eta,\zeta)$. 

\begin{figure}
\centeredepsfbox{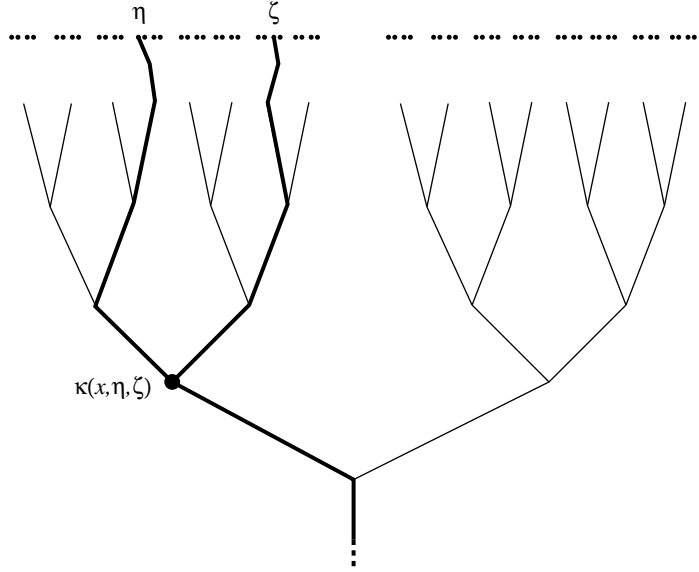}
\caption{The point $\kappa(x,\eta,\zeta)$ is the intersection of the three
sides of the triangle $\Delta(-\infinity,\eta,\zeta)$ in the tree $\tau \subset
X_n$ correponding to $x$.}
\label{FigureTreeBarycenter}
\end{figure}

\end{proof}

We now give a property of uniform biconvergence groups which is similar in
spirit to what holds for uniform convergence groups.  

\begin{proposition}[Contraction property]
\label{proposition:contraction}
Let $G$ be a finitely generated group quasi-isometric to $X_n$. Let
$(\rho_\low,\rho^u) \from G \to (\Bilip(\reals),\Bilip(\Q_n))$ be the uniform
biconvergence action given by Proposition \ref{proposition:biconvergence}. For
every compact subset $K\subset \Q_n$ and every open subset $U\subset \Q_n$,
there exists $g\in G$ such that $g \cdot K\subset U$.
\end{proposition}

\begin{proof} The appropriate generalization of this 
proposition should hold for
any uniform biconvergence group, but here we will make use of the fact
that $\rho^u$ is a uniform quasisimilarity action.

Choose $x \in \reals$ and $\zeta \in \Q_n$, and choose a one-to-one convergent
sequence $\zeta_1, \zeta_2, \ldots \to \zeta$ in $\Q_n$ such that $\zeta_i \ne
\zeta$ for all $i \in \naturals$. Choose a compact fundamental domain $K \subset
T(\reals,\Q_n,\Q_n)$ for the action of $G$. Choose $g_i \in G$ so that
$g_i^\inverse \cdot (x,\zeta_i,\zeta) \in K$. Pass to a subsequence so that
$g_i^\inverse \cdot (x,\zeta_i,\zeta)$ converges to $(y,\eta,\omega) \in K$. Of
course $\eta \ne \omega$, because $(y,\eta,\omega) \in T(\reals,\Q_n,\Q_n)$. Note
that $d(g_i^\inverse \cdot \zeta_i,g_i^\inverse\cdot \zeta) \to
d(\eta,\omega)\ne 0$, while $d(\zeta_i,\zeta) \to 0$. It follows that the
stretch interval $[a_i,b_i]$ for the action of $g_i$ on $\Q_n$ converges to
$0$, and so for sufficiently large $i$ the map $g_i$ contracts distance
uniformly in $\Q_n$, with a contraction factor $b_i$ that goes to $0$ as $i
\to\infinity$. On the other hand, we also know that $d(\omega,g_i^\inverse
\cdot\zeta)$ converges to $0$, and so 
\begin{align*}
d(g_i \cdot \omega,\zeta) &= d(g_i \cdot \omega, g_i^{\vphantom{\inverse}}
\cdot (g_i^\inverse\cdot\zeta)) \\ 
&\le b_i \, d(\omega,g_i^\inverse \cdot \zeta) \\
&\to 0 \quad\text{ as }\quad i \to \infinity
\end{align*} 
It follows that the action of the sequence $(g_i)$ on $\Q_n$ converges uniformly
on compact sets to the constant map with value $\zeta$. Choosing $\zeta$ to be
any point in $U$, and choosing $K \subset \Q_n$ to be any compact set, the
proposition follows.
\end{proof}

\begin{Remark} 
In the above proof, we know by Proposition \ref{prop:BilipAction} that the
stretch intervals for the action of $g_i$ on $\reals$ and on $\Q_n$ are
inversely related, and so the stretch interval for $g_i^\inverse$ on $\reals$
converges to $0$. Since $d(g_i^\inverse(x),y) \to 0$ it follows that
$g_i^\inverse$ converges uniformly on compact sets to the constant map with
value $y$. In some sense the sequence $(g_i)$ therefore has ``source-sink''
dynamics, with source $y \in \R$ and sink $\zeta \in \Q_n$. This gives another
analogy between convergence groups and biconvergence groups.
\end{Remark}

\subsection{Two applications of the biconvergence property}

Our first application shows that the representation $\theta \from G \to
\Aff(\reals)$ is virtually faithful.

\begin{proposition}[$\theta$ is virtually faithful] 
The homomorphisms $\rho_\ell \from G \to \Bilip(\reals)$ and $\theta \from G \to
\Aff(\reals)$ each have finite kernel.
\label{proposition:virtfaithful}
\end{proposition}

\begin{proof}
Since $\rho_\ell$ and $\theta$ are conjugate in $\Bilip(\reals)$, they
have the same kernel. So we need only prove that $\rho_\ell$ has finite kernel.

Pick a fixed clone $C \subset \Q_n$. Let $A = 0 \cross 1 \cross C$, a compact
subset of $T(\R,\R,\Q_n)$. Note that $C$ is also an open subset of $\Q_n$.

Suppose the proposition is false. Then there exists an arbitrarily large subset
$\{g_1,\ldots,g_k\} \subset G$ acting trivially on $\R$. 
We shall construct $g \in G$, depending on $g_1,\ldots,g_k$, such that $(g g_i
g^\inverse) \cdot A \intersect A \ne \emptyset$ in $T(\R,\R,\Q_n)$, for
$i=1,\ldots,k$. Since $k$ is arbitrarily large, this contradicts proper
discontinuity of the action of $G$ on $T(\R,\R,\Q_n)$ (Proposition
\ref{proposition:biconvergence}).

Pick any point $\zeta_0 \in \Q_n$. Let $\zeta_i = g_i \cdot \zeta_0$. Let $K$ be
the smallest clone containing $\zeta_0, \zeta_1, \ldots, \zeta_k$, in particular
$K$ is compact. Furthermore, $(g_i \cdot K) \intersect K \ne \emptyset$ for
$i=1,\ldots,k$.

By Proposition \ref{proposition:contraction} applied to the compact set $K$ and
the open set $C$, there exists $g \in G$ such that $g \cdot K \subset C$ in
$\Q_n$. Thus, $(g g_i g^\inverse) \cdot C \intersect C \ne \emptyset$ in $\Q_n$,
for each $i=1,\ldots,k$. Since $g g_i g^\inverse$ acts as the identity on
$\reals$ it follows that $(g g_i g^\inverse) \cdot A \intersect A \ne \emptyset$
in $T(\R,\R,\Q_n)$, contradicting proper discontinuity as indicated above. 
\end{proof}

Our next application of the biconvergence property shows that the stretch group
of $\theta(G) \subset \Aff(\reals)$ is infinite cyclic.

\begin{proposition}[stretch group is cyclic] 
Let $\Gamma = \theta(G)$. The group $\Stretch(\Gamma)$ is infinite cyclic.
\label{proposition:cyclic:stretch}
\end{proposition}

\begin{proof}
Recall that the representation $\theta \from G \to \Aff(\reals)$ is obtained
from $\rho_\ell \from G \to \Bilip(\reals)$ by a bilipschitz conjugacy: there
exists a bilipschitz map $\phi \from \reals \to \reals$ such that 
$$\theta(g) = \phi \composed \rho_\ell(g) \composed \phi^\inverse
$$
for all $g \in G$.

Let $s(g) = \Stretch(\theta(g))$. Let $[a_\ell(g),b_\ell(g)]$ be the stretch
interval of $\rho_\ell(g)$ acting on $\R$, and let $[a^u(g),b^u(g)]$ be the
stretch interval of $\rho^u(g)$ acting on $\Q_n$. We know by Proposition
\ref{prop:BilipAction} that there exists $L_1 \ge 1$ such that
$$[a^u(g),b^u(g)] \cdot [a_\ell(g),b_\ell(g)] \subset
\bigl[\frac{1}{L_1},L_1\bigr]
$$
We also know since the conjugating map $\phi$ is bilipschitz that
\begin{align}
\frac{[a_\ell(g),b_\ell(g)]}{s(g)} \subset \bigl[ \frac{1}{L_2},L_2 \bigr]
\tag*{$(*)$}
\end{align}
where $L_2$ is the square of the bilipschitz constant of $\phi$. Therefore,
setting $L_3 = L_1 L_2$ we have
\begin{align} [a^u(g),b^u(g)] \cdot s(g) \subset \bigl[ \frac{1}{L_3},L_3
\bigr] \tag*{$(**)$}
\end{align}

Suppose the proposition is false. Then there exists a sequence $g_i \in G$ such
that $s(g_i) \ne 1$ but $s(g_i) \to 1$ as $i \to \infinity$. Since $s(g_i) \ne
1$ it follows that $\theta(g_i)$ has a unique fixed point in $\reals$, and so
$\rho_\ell(g_i)$ has a unique fixed point in $\reals$.

We claim also that $\rho^u(g_i)$ has a unique fixed point in $\Q_n$. To see
this, since $s(g_i) \ne 1$ we may choose a sufficiently high power $g_i^k$ so
that $s(g_i^k) \not\in \bigl[ \frac{1}{L_3},L_3 \bigr]$. It follows that
$$1 \not\in [a^u(g_i^k),b^u(g_i^k)]
$$
and therefore either $\rho^u(g_i^k)$ or $\rho^u(g_i^{-k})$ is a contraction
mapping of $\Q_n$ and so has a unique fixed point, and so $\rho^u(g_i)$
has a unique fixed point.

Let $x_i \in \reals$ be the unique fixed point of $\rho_\ell(g_i)$, and let
$\zeta_i \in \Q_n$ be the unique fixed point of $\rho^u(g_i)$. For the rest of the
proof we use only the actions $\rho_\ell \from G \to \Bilip(\reals)$ and
$\rho^u \from G \to \Bilip(\Q_n)$, and so we may unambiguously use the ``dot''
notation for these group actions.

Since the action of $G$ on $T(\R,\R,\Q_n)$ is cocompact, we may pick a compact
fundamental domain~$A$.

Pick any point $w \in \reals$ distinct from all the points $x_i$. Choose
$h_i \in G$ so that $h_i(x_i,w,\zeta_i) \in A$. Replacing $g_i$ by its conjugate
$h_i g_i h_i^\inverse$, and replacing $x_i,w,\zeta_i$ by their images under
$h_i$, we may assume that $(x_i,w_i,\zeta_i) \in A$ for some sequence
$w_i$; since $s$ is invariant under conjugation we still have $s(g_i) \to 1$. By
compactness of~$A$ we may pass to a subsequence so that $x_i \to x$ in $\reals$
and $\zeta_i \to \zeta$ in $\Q_n$. 

Since $s(g_i) \to 1$, using $(*)$ and $(**)$ it follows that for all
sufficiently large~$i$ we have
$$[a_\ell(g_i),b_\ell(g_i)] \subset \bigl[\frac{1}{L}, L \bigr]
$$
and
$$[a^u(g_i),b^u(g_i)] \subset \bigl[ \frac{1}{L}, L \bigr]
$$
where $L = 2\,\Max(L_2, L_3)$.

Choose a point $y \ne x$ in $\reals$. We will contradict proper
discontinuity of the action of $G$ on $T(\R,\R,\Q_n)$ by showing that
$(g_i \cdot x,g_i \cdot y,g_i \cdot \zeta)$ stays in a compact subset of
$T(\R,\R,\Q_n)$. 

We know that $g_i \cdot x \to x$ in $\reals$, because 
\begin{align*}
d(x,g_i \cdot x) &\le d(x,x_i) + d(x_i, g_i \cdot x) \\
                 &= d(x,x_i) + d(g_i \cdot x_i, g_i \cdot x) \\
                 &\in d(x,x_i) + d(x,x_i) \cdot \bigl[\frac{1}{L}, L \bigr] \\
                 &= d(x,x_i) \bigl[1 + \frac{1}{L}, 1 + L \bigr] \\
                 &\to 0 \quad\text{as}\quad i \to \infinity
\end{align*} 
because $d(x,x_i) \to 0$. The same argument shows that $g_i \cdot \zeta \to
\zeta$ in $\Q_n$.

To complete the proof, we will show that $g_i \cdot y$ stays in a compact subset
of $\reals$ disjoint from $x$, from which it follows that
$(g_i \cdot x,g_i \cdot y,g_i \cdot \zeta)$ stays in a compact subset of triple
space, a contradiction as noted above.

We have
\begin{align*}
d(x_i,g_i \cdot y) &= d(g_i \cdot x_i,g_i \cdot y) \\
&\subset \bigl[ \frac{1}{L},L \bigr] \cdot d(x_i,y) \\
&\subset \bigl[ \frac{1}{2L},2L \bigr] \cdot d(x,y)
\end{align*}
for sufficiently large $i$, because of the fact that $d(x_i,y) \to
d(x,y)$. Since $d(x_i,x) \to 0$ it follows that
$$ d(x,g_i \cdot y) \in \bigl[ \frac{1}{3L},3L \bigr] \cdot d(x,y)
$$
for sufficiently large $i$. The set of points in $\reals$ whose distance from
$x$ is in the interval $[1/3L,3L] \cdot d(x,y)$ is a compact subset of
$\reals$ disjoint from $x$, and this compact subset contains $g_i \cdot y$
for all sufficiently large $i$, as desired.
\end{proof}

\begin{Remark} The inequality $(*)$ is \emph{precisely} where we need Theorem
\ref{theorem:qshink}, which gives a bilipschitz conjugacy between
$\rho_\ell$ and $\theta$, instead of the quasisymmetric conjugacy provided
by Hinkkanen's original theorem.
\end{Remark}

\vfill\break

\section{Finishing the proof of Theorem A}
\label{section:proof}

Suppose $G$ is quasi-isometric to $\BS(1,n)$. By Theorem
\ref{theorem:homomorphism} we have a quotient group $\Gamma = G / N$, where $N$
is a finite normal subgroup of $G$, and we have a diagram of short exact
sequences 
$$
\xymatrix{
1 \ar[r] & A=A(\Gamma) \ar@{^{(}->}[r] \ar@{^{(}->}[d] & \Gamma \ar[r]
\ar@{^{(}->}[d] & \Stretch(\Gamma)=\langle t \rangle \ar[r] \ar@{^{(}->}[d] & 1 \\
1 \ar[r] & \Isom(\R) \ar@{^{(}->}[r] & \Aff(\R) \ar[r] & \Stretch(\R) \ar[r] & 1
}
$$
%
%
where $\left< t \right>$ is infinite cyclic. Choosing a splitting, we may
regard $t$ as an element of $\Aff(\reals)$, and replacing $t$ with $t^\inverse$
if necessary we may assume $\Stretch(t) > 1$. 

We have the following additional properties of $\Gamma$:
\begin{itemize}
\item $\Gamma$ is solvable, because $\Aff(\reals)$ is solvable.
\item $\Gamma$ has a torsion free subgroup of index at most $2$, namely
$$\Gamma_+ = \Gamma \intersect \Aff_+(\reals)$$ 
\item $\Gamma$ is quasi-isometric to $\BS(1,n)$.
\end{itemize}
Since finite presentability is a quasi-isometry invariant (\cite{GhysHarpe:afterGromov},
Proposition 10.18) we also have: 
\begin{itemize}
\item $\Gamma$ is finitely presented.
\end{itemize}
Note that there are finitely generated subgroups of $\Aff(\reals)$ which are not 
finitely presented; once such group is described in \cite{Strebel:fpsolvable}.

We may now quote the following theorem of Bieri-Strebel, taken from
\cite{Strebel:fpsolvable} (see also \cite{BieriStrebel:almostfpsolvable}):

\begin{theorem}[Bieri-Strebel]
\label{theorem:bieri}
Let $\Gamma$ be a finitely presented solvable group, and suppose that $\Gamma$
has an HNN presentation of the form
$$\Gamma = A*_\Z = \left<A,t \suchthat t A_1 t^\inverse = A_2 \right>
$$
where $A_1, A_2$ are subgroups of $A$. Then there is another HNN presentation
$$\Gamma = B*_\Z = \left<B,t \suchthat t B t^\inverse = B'\right>
$$
where $B$ is a finitely generated subgroup of $A$ and $B'$ is a subgroup of $B$.
\end{theorem}

For completeness, here is a quick proof suggested to us by T. Delzant. 

\begin{proof} Let $K$ be a finite complex with fundamental group $\Gamma$ and
universal covering map $p \from \widetilde K \to K$. Let $T$ be a Bass-Serre tree
for the given HNN decomposition of $\Gamma$. Let $f \from \widetilde K \to T$ be
a $\Gamma$-equivariant map, transverse to the midpoint of each edge of $T$.
Pulling back those midpoints gives a $\Gamma$-equivariant 1-complex $\widetilde
L \subset \widetilde K$. Let $L = p(\widetilde L)$, a 1-complex in $K$.
Construct a new $\Gamma$-tree $T'$ whose edges correspond to components of
$\widetilde L$ and whose vertices correspond to components of $\widetilde K -
\widetilde L$. The vertex and edge stabilizers of $T'$ are subgroups of vertex
and edge stabilizers of $T$. Moreover, if $v$ is a vertex of $T'$ corresponding
to a component $\widetilde C$ of $\widetilde K - \widetilde L$, and if $C =
p(\widetilde C)$, then $\Stab(v)$ is isomorphic to the image of the inclusion
induced map $\pi_1(C) \to \pi_1(K)$, and so $\Stab(v)$ is finitely generated. 
Passing to the quotient graph of groups obtained from the action of $\Gamma$ on
$T'$, and collapsing subgraphs in the appropriate manner, one obtains an HNN
decomposition $\Gamma = \left<B,t' \suchthat t' B_1 t'{}^\inverse = B_2 \right>$
where $B$ is a finitely generated subgroup of $A$, $B_1, B_2 \subgroup B$, and
$t'$ is conjugate into the infinite cyclic subgroup $\left< t \right>$. If both
$B_1$ and $B_2$ are proper subgroups of $B$, then $\Gamma$ has a free subgroup of
rank $\ge 2$; this follows from the normal forms theorem for HNN decompositions,
or from a ping-pong argument on the Bass-Serre tree. Since $\Gamma$ is solvable
it has no free subgroups of rank $\ge 2$, and so one of $B_1, B_2$ equals $B$.
\end{proof} 

Applying the above theorem to $\Gamma$, we have a finitely generated subgroup
$B \subset A(\Gamma)$, and an HNN decomposition
$$\Gamma = \left< B, t \suchthat t b t^\inverse = \phi(b), \, \forall b \in B
\right>
$$
where $\phi \from B \to B$ is an injective endomorphism, and $B$ is a finitely
generated subgroup of $\Isom(\reals)$.

Recall that $\Gamma_+ = \Gamma \intersect \Aff_+(\reals)$, and let $B_+ = B
\intersect \Isom_+(\reals)$. The group $B_+$ is a finitely generated subgroup
of $\Isom_+(\reals) = \Transl(\reals)$, and so $B_+$ is free abelian of some
finite rank $r \ge 1$. The indices $[\Gamma : \Gamma_+]$, $[B : B_+]$ are both
$\le 2$, and so we fall into three cases:
\begin{description}
\item[Case 1.] $[\Gamma : \Gamma_+] = [B : B_+] = 1$.
\item[Case 2.] $[\Gamma : \Gamma_+] = 2$ and $[B : B_+] = 1$.
\item[Case 3.] $[\Gamma : \Gamma_+] = [B : B_+] = 2$.
\end{description}
In case (1) we obviously have $t \in \Aff_+(\reals)$, that is, $t$ preserves
orientation of $\reals$. In case (3), if $t \in \Aff(\reals) - \Aff_+(\reals)$,
in other words if $t$ reverses orientation of $\reals$, then we can replace $t$
by its product with an element of $B - B_+$, a reflection of $\reals$; and hence
we may assume in case (3) that $t\in\Aff_+(\reals)$. In case (2) we necessarily
have $t \in \Aff(\reals) - \Aff_+(\reals)$. Setting $s=t^2$ and $\psi = \phi^2$
in case (2), or $s=t$ and $\psi=\phi$ in cases (1) and (3), we have an HNN
decomposition
$$\Gamma_+ = \left< B_+, s \suchthat s b s^\inverse = \psi(b), \, \forall b \in
B_+ \right>
$$
where $s \in \Aff_+(\reals)$. Since $B_+ \homeo \integers^r$
and $\psi(B_+) \homeo B_+$, it follows that $\psi(B_+)$ has finite index in
$B_+$. Let $I = [B_+ : \psi(B_+)]$. 

Next we extract some homological information about $\Gamma_+$ (see
\cite{Brown:cohomology}). Recall that if $K$ is a finitely generated, virtually
torsion free group, the \emph{virtual cohomological dimension} of $K$, denoted
$\vcd(K)$, is the cohomological dimension of any finite index, torsion free
subgroup $K' \subgroup K$; by Serre's Theorem (\cite{Brown:cohomology} Theorem
VIII.3.1) this number is independent of the choice of $K'$. 

\begin{lemma} For any injective homomorphism $\psi \from \Z^r \to \Z^r$, the
HNN group $K = \left< \Z^r,t \suchthat t b t^\inverse = \psi(b), \, \forall b \in
\Z^r \right>$ has virtual cohomological dimension equal to $r+1$.
\label{lemma:vcd}
\end{lemma}

\begin{proof} We first construct a compact Eilenberg-Maclane space for $K$ by
mimicking the construction in \cite{FarbMosher:BSOne} of an
Eilenberg-Maclane space for $\BS(1,n)$.

Extend the endomorphism $\psi \from \Z^r \to \Z^r$ to a linear isomorphism
$\psi \from \reals^r \to \reals^r$ which commutes with the action of $\Z^r$. Let
$T^r = \R^r / \Z^r$ be the $r$-dimensional torus. The map $\psi$ descends to a
self covering map $\Psi \from T^r \to T^r$. Let $I$ be the index
$[\Z^r : \psi(\Z^r)]$, and so $I$ is the determinant of $\psi$ acting on
$\reals^r$, and $I$ is the degree of the covering map $\Psi$. Let $C$ be the
mapping torus of $\Psi$,
$$C = T^r \cross [0,1] \bigm/ (x,0) \sim (\Psi(x),1)
$$
Note that $\pi_1(C) \homeo K$. Let $X$ be the universal covering space of $C$.
Let $T_I$ be the tree constructed in \S\ref{section:background}, a homogeneous
directed tree with one incoming edge and $I$ outgoing edges at each vertex. By
construction, the space $X$ fibers over $T_I$ with fiber $\reals^r$. Since
$T_I$ is contractible we have a homeomorphism
$$X \homeo \reals^r \cross T_I
$$
In particular $X$ is contractible, and so $C$ is a compact Eilenberg-Maclane
space for the group $K$.

Since $K$ is torsion free we have $\vcd(K)=\cd(K)$. By \cite{Brown:cohomology}
corollary VIII.7.6, the number $\cd(K)$ is the maximum dimension $k$ such that
$H^k_c(X;\integers) \ne 0$, where $H^k_c$ denotes cohomology with compact
supports. By \cite{Spanier} (p. 341, corollary 9 and p. 360, exercise 6), there
is a K\"unneth formula
$$H^k_c(X;\integers) \homeo \Sum_{i=0}^k H^i_c(\reals^r;\integers) \tensor
H^{k-i}_c(T_I;\integers)
$$
from which it follows that
$$H^k_c(X;\integers) \homeo \begin{cases}
0 & \text{if $k \ne r+1$}\\
H^1_c(T_I;\integers) & \text{if $k = r+1$}
\end{cases}
$$
Since $T_I$ has infinitely many ends, it follows that $H^1_c(T_I;\integers)$ is
an infinite rank abelian group, and so $\vcd(K) = r+1$.
\end{proof}

\begin{Remark}
This proof gives extra information about the group $K$, namely a complete
computation of the cohomology groups $H^k_c(X;\integers) \homeo H^k(K;\Z K)$,
which by \cite{Gersten:dimension} are quasi-isometry invariants of $K$. Note in
particular that $K$ is a duality group in the sense of Bieri and Eckmann (see
\cite{Brown:cohomology} \S VIII), which means that, setting $D = H^{r+1}(K;\Z K)
\homeo H^{r+1}_c(X;\Z)$, we have
$$H^i(K,M) \homeo H_{r+1-i}(K,D \tensor M) \text{ for all
$K$-modules $M$.}
$$
\end{Remark}
\medskip

Now we complete the proof of Theorem A.

Applying Lemma \ref{lemma:vcd} 
we have the following computations of virtual
cohomological dimensions:
$$\vcd(\Gamma_+) = r+1 \text{ and } \vcd(\BS(1,n)) = 2
$$
By a theorem of Gersten \cite{Gersten:dimension} and Block-Weinberger \cite{BlockWeinberger}, if two
virtually torsion free groups have finite Eilenberg-Maclane spaces, and if the
groups are quasi-isometric, then they have the same $\vcd$. It follows that
$\vcd(\Gamma_+) = \vcd(\BS(1,n))$, and so $r = 1$ and $B_+$ is an infinite cyclic
group $B_+ = \left<b\right>$. The endomorphism $\psi \from \integers \to
\integers$ preserves orientation, and so $\psi(b)=b^m$ for some $m \ge 1$. If
$m=1$ then $\Gamma_+ \homeo \integers^2$, but $\integers^2$ is not
quasi-isometric to $\BS(1,n)$---for example, $\integers^2$ has quadratic growth
but $\BS(1,n)$ has exponential growth. It follows that $m \ge 2$, and $\Gamma_+$
is isomorphic to $\BS(1,m)$. Applying the main result of \cite{FarbMosher:BSOne} it follows
that $\Gamma_+$ is abstractly commensurable to $\BS(1,n)$, in fact $m,n$ are
positive powers of the same positive integer.

As for the group $\Gamma$ itself, we have already shown that $\Gamma$ has a
subgroup $\Gamma_+$ isomorphic to $\BS(1,m)$ of index $\le 2$, and $\BS(1,m)$ is
abstractly commensurable to $\BS(1,n)$; therefore $\Gamma$ is abstractly
commensurable with $\BS(1,n)$. This completes the proof of Theorem
A. \qed

\bigskip

Here is a corollary:

\begin{corollary} Let $G$ be a finitely generated, torsion free group which is
quasi-isometric to $\BS(1,n)$ for some $n \ge 2$. Then $G \homeo \BS(1,k)$ for
some integer $k$ with $\absvalue{k} \ge 2$, such that $\BS(1,n)$ and $\BS(1,k)$
are abstractly commensurable.
\label{corollary:torsionfree}
\end{corollary}

\begin{proof} Apply Theorem A to get a short exact sequence 
$$1 \to N \to G \to \Gamma \to 1
$$ 
with $N$ finite. Since $G$ is torsion free it follows that $N = \{\Id\}$ and $G
\homeo \Gamma$. In the proof of Theorem A, it also follows that $B=B_+$ is
infinite cyclic, $B = \left< b \right>$. Thus we fall into case 1 or 2 in the
proof of Theorem A. The injective endomorphism $\phi\from B\to B$ must have the
form $\phi(b)=b^k$ for some integer $k \ne 0$. If $k=\pm 1$ then $\Gamma$ has
polynomial growth, a contradiction as above. Therefore, $\absvalue{k} \ge 2$ and
$G \homeo \BS(1,k)$. Since $\BS(1,k)$ contains $\BS(1,k^2)$ with index $2$, and
since $k^2 \ge 2$, the main result of \cite{FarbMosher:BSOne} shows that
$\BS(1,k^2)$ is abstractly commensurable to $\BS(1,n)$, and so $G \homeo
\BS(1,k)$ is abstractly commensurable to $\BS(1,n)$.
\end{proof}

In the proof of Theorem A, recall the three cases in the analysis of the
group $\Gamma$. As the proof of the corollary shows, in case 1 we have $\Gamma
\homeo \BS(1,k)$ for some $k \ge 2$, and in case 2 we have $\Gamma \homeo
\BS(1,k)$ for some $k \le -2$. 

In case 3, $B$ is the infinite dihedral group $B = \left<a,r \suchthat r^2
= 1, rar = a^\inverse\right>$. To enumerate the possibilities for $\Gamma$, it
suffices to enumerate the injective, nonsurjective, orientation preserving
endomorphisms $\phi \from B \to B$ (orientation preserving means that $\phi$
fixes each of the two ends of the group $B$). Fix the representation $B
\to \Aff(\reals)$ given by $a \cdot x = x+2$, $r \cdot x = -x$. The reflections
of $B$ may be enumerated as $r_i = a^{-i}r$, $i \in \integers$, and $r_i$ is
the reflection about the point $i \in \reals$. Up to conjugation by an
automorphism of $B$, there are two infinite sequences of endomorphisms $\phi$,
depending on whether $\phi$ fixes some reflection, as follows.
\begin{description}
\item[Case 3.i.] If $\phi$ fixes a reflection then, replacing $r$
by the fixed reflection, there exists an integer $m \ge 2$ such that
$\phi(a)=a^m$, $\phi(r)=r$. We may represent $\phi$ in $\Aff_+(\reals)$ as the
expansion with fixed point $0$ and stretch $m$, that is, $\phi(x)=mx$. We
therefore have
$$\Gamma = B*_\phi = \left< a,r,t \suchthat r^2=1, rar^\inverse = a^\inverse,
tat^\inverse = a^m, trt^\inverse = r \right>
$$
\item[Case 3.ii.] If $\phi$ does not fix a reflection then, up to a replacement
of
$r$, there is an odd integer $m=2k+1 \ge 3$ such that $\phi(a)=a^m$ and $\phi(r)
= a^{-k} r$. We may represent $\phi$ in $\Aff_+(\reals)$ as the expansion with
fixed point $-1/2$ and stretch $m$, that is $\phi(x) = m(x+1/2)-1/2$. It is
clear from this description that $\phi$ indeed fixes no reflections; one can also
compute directly that $\phi(r_i) = r_{2ki+k+i}$, and clearly there is no $i \in
\integers$ such that $2ki+k+i=i$. We have:
$$\Gamma = B*_\phi = \left< a,r,t \suchthat r^2=1, rar^\inverse = a^\inverse,
tat^\inverse = a^m, trt^\inverse = a^{-k}r \right>
$$
\end{description}

\section{Final comments: The Sullivan-Tukia Theorem}
\label{section:questions}

The Sullivan-Tukia Theorem says that a uniformly quasiconformal
subgroup of $\QC(S^2) \homeo \QI(\hyp^3)$ is quasiconformally conjugate
into the conformal group $\Conf(S^2) \homeo \Isom(\hyp^3)$.

Here is one possible way to formulate an analogue of this theorem for
$\BS(1,n)$. The analogue of the group $\QC(S^2) \homeo \QI(\hyp^3)$ 
is $\QI(\BS(1,n)) \homeo\QI(X_n)
\homeo \Bilip(\reals) \cross \Bilip(\Q_n)$. The analogue of a uniformly
quasiconformal subgroup of $\QC(S^2)$ is a ``bounded'' subgroup $G \subgroup
\QI(X_n)$, which means that for some $K \ge 1$, each $g \in G$ is
represented by a $(K,C)$-quasi-isometry of $X_n$ (where $C$ may depend on $g$).
Equivalently, there is a constant $K' \ge 1$ such that for each $g \in G$,
letting $\rho_\ell(g) \in \Bilip(\reals)$, $\rho^u(g) \in \Bilip(\Q_n)$ be the
projections, we have 
$$[a(\rho_\ell(g)),b(\rho_\ell(g))] \cdot [a(\rho^u(g)),b(\rho^u(g))] \subset
[1/K',K']
$$
What is the analogue of $\Conf(S^2)$? One guess is $\Isom(X_n)$. However,
the Sullivan-Tukia Theorem is false in this case, because 
if $n=m^2$ then $\BS(1,m)$ is a
bounded subgroup of $\QI(X_n)$, but there does not exist an isometric
action of $\BS(1,m)$ on~$X_{n}$. 

A more reaonable guess is that if $n\ge 2$ is an integer which is not a
proper power, then every bounded subgroup of $\QI(X_n)$ is conjugate to a
subgroup of $\Isom(X_n)$. One possible approach is to prove a version of
Theorem \ref{theorem:qshink} for 
$\Q_n$, namely that any uniform quasisimilarity
group action on $\Q_n$ is bilipschitz conjugate to a similarity action; for
this to be true it is necessary that $n$ not be a proper power. However, we do
not know how to prove this.

\vfill\break

\bibliographystyle{amsalpha}
\bibliography{mosher}

\bigskip

\noindent
Benson Farb:\\
Department of Mathematics, Univeristy of Chicago\\
5734 University Ave.\\
Chicago, Il 60637\\
E-mail: farb@math.uchicago.edu
\medskip

\noindent
Lee Mosher:\\
Department of Mathematics, Rutgers University, Newark\\
Newark, NJ 07102\\
E-mail: mosher@andromeda.rutgers.edu

\end{document}